\documentclass[10pt]{article}
\usepackage{wrapfig}
\usepackage{graphicx} 
\usepackage{amsmath, amssymb}
\usepackage[english]{babel} 
\usepackage{amsthm} 
\usepackage{amssymb,amsmath}
\usepackage{subcaption}
\usepackage{bbm} 
\usepackage{enumerate}
\usepackage{graphicx}
\usepackage{float}
\usepackage{refcount} 
\usepackage{enumitem} 
\usepackage{mathtools} 
\usepackage[dvipsnames]{xcolor} 

\usepackage{soul} 
\sethlcolor{lightgray} 
\usepackage[final, linkcolor=blue]{hyperref}
\hypersetup{
	colorlinks   = true, 
	urlcolor     = blue, 
	linkcolor    = blue, 
	citecolor   = red 
}
\usepackage{mathrsfs} 
\usepackage{cite} 
\usepackage[a4paper]{geometry}
\usepackage{multicol}
\usepackage[bottom]{footmisc}

\usepackage{url}

\usepackage{tikz-cd}
\usepackage{amsmath}
\usepackage{tikz}
\usepackage{pgfplots}
\usetikzlibrary{positioning}
\tikzstyle{line} = [draw, -latex']
\pgfplotsset{compat=newest}
\pgfplotsset{plot coordinates/math parser=false}
\newlength\figureheight
\newlength\figurewidth

\newcommand{\Ndata}{N_{\text{batch}}}
\newcommand{\Nbatch}{N_{\text{batch}}}

\newcommand{\N}{\mathbb{N}}

\newcommand{\R}{\mathbb{R}}

\newcommand{\D}{\mathcal{D}}

\newcommand\numberthis{\addtocounter{equation}{1}\tag{\theequation}} 

\renewcommand{\epsilon}{\varepsilon}
\newcommand{\eps}{\epsilon}

\makeatletter
\def\blfootnote{\gdef\@thefnmark{}\@footnotetext}
\makeatother

\newtheorem{theorem}{Theorem}

\theoremstyle{definition}

\newtheorem{remark}[theorem]{Remark}

\newcounter{subexample}
\renewcommand{\thesubexample}{\theexample.\arabic{subexample}}

\newenvironment{subexample}[1][]{%
	\refstepcounter{subexample}%
	\par\noindent\textbf{Example \thesubexample}%
	\ifx\relax#1\relax\else\ \textit{(#1)}\fi.
}{\par}

\newcommand{\setexample}[1]{\setcounter{example}{#1}\setcounter{subexample}{0}}

\numberwithin{equation}{section} 

\makeatletter
\let\@fnsymbol\@arabic
\makeatother

\title{Tracking Finite-Time Lyapunov Exponents to Robustify Neural ODEs}

\author{
    Christian Kuehn\thanks{TU Munich, School of CIT, Department of Mathematics, Boltzmannstra\ss e 3, 85748 Garching bei M\"unchen. E-mail: \texttt{christian.kuehn@tum.de}.} 
    \and 
    Tobias W\"ohrer\thanks{TU Wien, Department of Mathematics, Institute of Analysis and Scientific Computing, Wiedner Hauptstra\ss e 8--10. E-mail: \texttt{tobias.woehrer@tuwien.ac.at}.}
}

\date{\today}

\begin{document}

\maketitle
\tableofcontents
\abstract{We investigate finite-time Lyapunov exponents (FTLEs), a measure for exponential separation of input perturbations, of deep neural networks within the framework of continuous-depth  neural ODEs. We demonstrate that FTLEs are powerful organizers for input-output dynamics, allowing for better interpretability and the comparison of distinct model architectures. We establish a direct connection between Lyapunov exponents and adversarial vulnerability, and propose a novel training algorithm that improves robustness by FTLE regularization. The key idea is to suppress exponents far from zero in the early stage of the input dynamics. This approach enhances robustness and reduces computational cost compared to full-interval regularization, as it avoids a full ``double'' backpropagation.
}

\section{Introduction}\label{sec:intro}
A fundamental property of nonlinear dynamics is the geometric deformation of phase space over time. By locally linearizing the flow along trajectories one obtains the variational equation. It measures how infinitesimal displacement vectors expand or contract along specific directions. Combining the original dynamics and the variational equation allows for the calculation of Lyapunov exponents (LEs). The LEs describe the magnitude of exponential separation and provide a precise measure of the system's sensitivity to perturbations.
LEs are a powerful tool for uncovering hidden structures in complex dynamical systems \cite{oseledec68,eckmann1986,DoanKarraschNguyenSiegmund}. Their finite-time counterparts (FTLEs) are widely used in non-autonomous settings, such as complex fluid dynamics and mixing processes where a full rigorous analysis is not feasible. They are used for the numerical identification of flow-invariant regions separated by Lagrangian coherent structures \cite{haller00,haller02,ShaddenLekienMarsden,haller11,NolanSerraRoss} and can be determined efficiently from data~\cite{WolfSwiftSwinneyVastano,GeistParlitzLauterborn,RosensteinCollinsDeLuca,KantzSchreiber}. The existence of a positive LE provides a quantitative indicator of chaotic behavior~\cite{Pesin,Young1,EckmannRuelle1,AlligoodSauerYorke}, and this idea can be extended to transient chaos using FTLEs~\cite{Abarbaneletal,LaiTel}. More recently, LEs have found initial applications in deep learning contexts. In this setting, Lyapunov exponents help to assess trainability \cite{engelken23,engelken2020lyapunov}, detect chaotic regimes \cite{storm2023finite,vockmeisel25} and link chaotic behavior to adversarial vulnerability \cite{pedraza2022lyapunov}. Recent years have seen deep neural networks, such as Residual Neural Networks (ResNets) \cite{resnet2016}, achieve state-of-the-art results in image recognition, natural language processing, and scientific discovery \cite{wang2023scientific}. A key theoretical advance has been the interpretation of ResNets as Euler discretizations of neural ordinary differential equations (nODEs). This continuous-depth perspective \cite{haber_multiscale2017, weinan2017} is one possible path to bridge the gap between discrete deep learning architectures and the rich theory of continuous-time dynamical systems.

To illustrate the connection to Lagrangian coherent structures, let us consider a nODE successfully trained to solve a binary classification problem. Then, the dynamics generate two flow-invariant regions corresponding to the two classes. In analogy to fluid dynamics, the classification dynamics could be understood as an \emph{unmixing process}. The invariant regions are separated by an ensemble of trajectories of codimension one, corresponding to ridges of locally maximal FTLE values, which ultimately determines the performance of the model. It is therefore of great interest to understand how the nODE generated input-output dynamics form the ridge between the classes. This viewpoint of ridges of FTLEs as organizing structures have already permeated several areas of nonlinear science. We build upon these insights, recognizing their power within the context of nODEs. 

In this work we investigate nODEs that solve two-dimensional classification tasks as a prototypical deep neural network application. The first goal is to track the FTLEs which show improved model interpretability by revealing the hidden invariant structures of the input-to-output process. The second goal is to use the acquired insights to improve the robustness of deep neural networks. 

\begin{figure}
	\centering 
    \begin{subfigure}{0.255\textwidth}	\centering 
		\includegraphics[width = \textwidth, trim={0 0.em 0 0},clip]{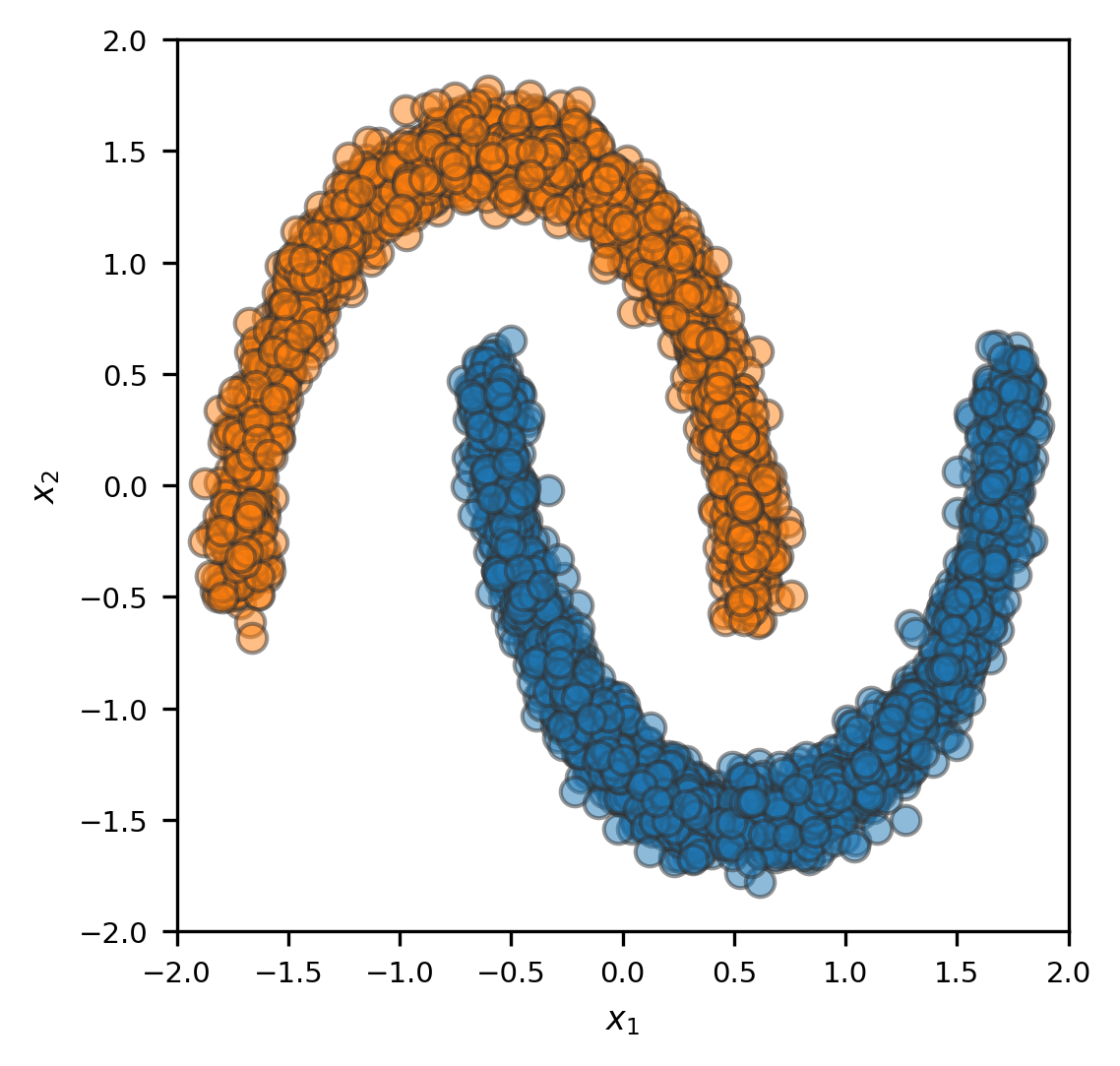}	
		\caption{\small Training data \label{subfig:trainingset}}
	\end{subfigure}
	\begin{subfigure}{0.31\textwidth}	\centering 
		\includegraphics[width = \textwidth, trim={0 0.em 0 1.6em},clip]{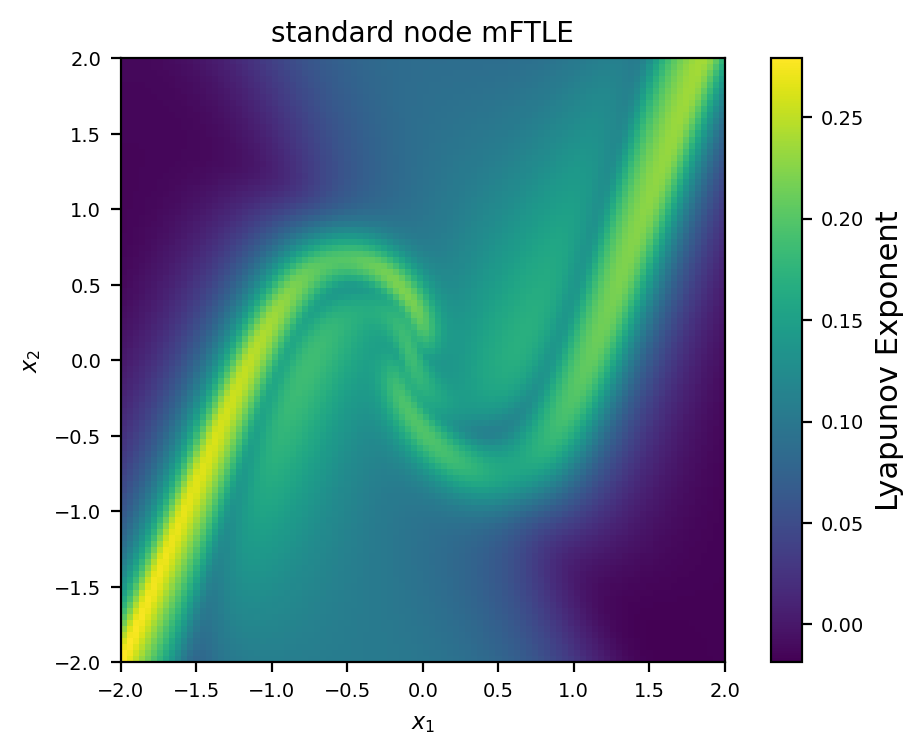}	
		\caption{\small Standard training \label{subfig:intronoreg}}
	\end{subfigure}
	\begin{subfigure}{0.31\textwidth}	\centering 
		\includegraphics[width = \textwidth, trim={0 0.em 0 1.6em},clip]{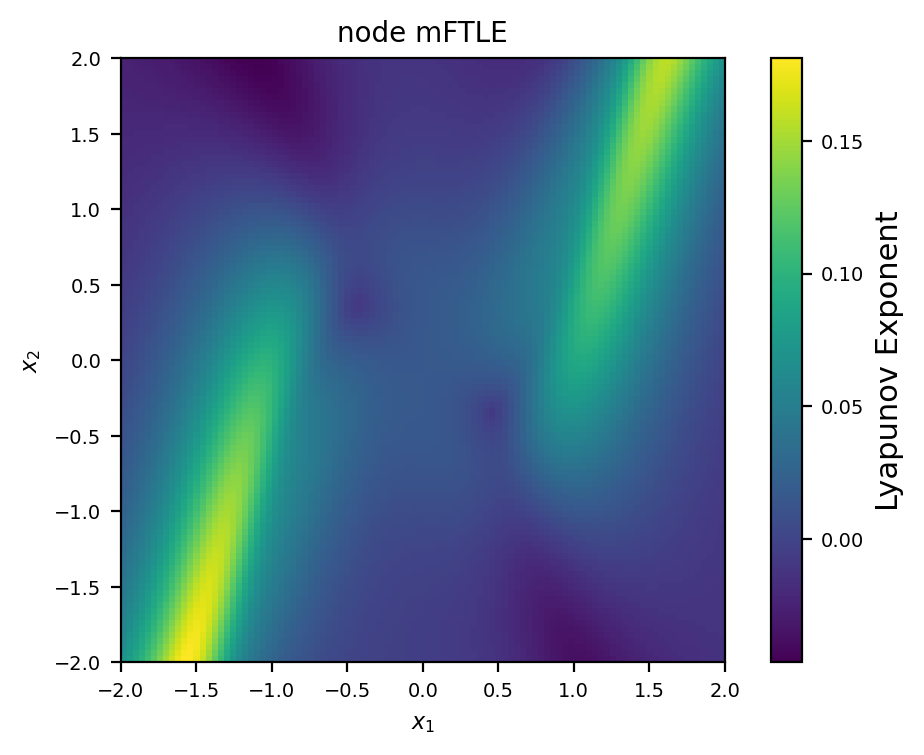}
		\caption{\small FTLE suppression\label{subfig:introreg}}
	\end{subfigure}
	\caption{\small Comparison of FTLEs of nODEs trained without and with active FTLE suppression. Figure~\ref{subfig:intronoreg} shows large FTLE ridge at the class boundary. The modified training of Figure~\ref{subfig:introreg} achieves a suppression of the FTLE ridge, reducing input sensitivity at the decision boundary.\label{fig:intro}}
\end{figure}

Within the framework of nODEs \cite{CRBD18, weinan2017}, we cover a broad class of model variations. When nODEs solve classification tasks, such as the moons dataset of Figure~\ref{subfig:trainingset}, their FTLE ridges reveal the geometry of decision boundaries and the underlying separation process, see Figure~\ref{subfig:intronoreg}.
Since Lyapunov exponents are by definition linked to input sensitivities,
they provide a powerful quantitative measure of adversarial vulnerability. Building on this understanding, we introduce a novel regularization term to modify the loss function that suppresses large Lyapunov exponents, see Figure~\ref{subfig:introreg},
leading to deep neural networks with enhanced adversarial robustness.

\subsection{General setting}

We consider \emph{neural ODEs} (nODEs) \cite{CRBD18} of the general form
\begin{equation}\label{eq:node}
    \begin{cases*}
    	\dot x(t) = f(\theta(t), x(t)),\quad t\in [0,T), \\
    	x(0) = x_0 \in \R^d,
    \end{cases*}
\end{equation}
where $T>0$, and the right-hand side includes explicit time-dependent parameters $\theta(t)$. For the scope of this paper (as is typical in many implementations), we assume the parameter function to be piecewise-constant in time, i.e.
\begin{equation}\label{eq:theta}
	\theta(t) = \sum_{k = 1}^K \theta_k \mathbbm{1}_{[\alpha_k,\beta_k)}(t),
\end{equation}
where $\bigcup_{k=1}^K [\alpha_k, \beta_k) = [0,T)$ and where $\theta_k$ are real-valued parameter collections specified below.
We refer to $K$ as the \emph{number of parameter layers}  of the nODE to emphasize the close connection to Residual Neural Network (ResNet) layers \cite{resnet2016, weinan2017, haber_multiscale2017}, see Remark~\ref{rem:resnetnode} below. The class of vector fields $f(\theta(t),\cdot):\R^d \to \R^d$ for which we consider solutions to \eqref{eq:node}, is given as
\begin{align}\label{eq:generalf}\begin{aligned}
	f(\theta(t),x) &= f_\ell(\theta(t),\cdot)\circ \cdots \circ f_1(\theta(t),x) \quad \text{with}\\
		f_i(\theta(t),y) &= V_i(t)\sigma_i\Big(W_i(t)y + b_i(t)\Big) + a_i(t),\quad i = 1, \ldots, \ell,
		\end{aligned}
\end{align}
where $\ell\in \N$ is the \emph{number of vector field layers}. The function $\theta(t)$ collects all the parameters for all parameter layers and vector field layers. Specifically, the parameter dependence of $f_i(\theta(t),\cdot): \R^{d_{i-1}} \to \R^{d_i}$ is given via the matrices
 $W_i(t) \in \R^{d_{i}^\text{hid}\times d_{i-1}}, V_i(t) \in \R^{d_i \times d_i^{\text{hid}}}$ and vectors $a_i(t) \in \R^{d_i}, b_i(t) \in \R^{d_i^\text{hid}}$, denoting $d_0:=d$, $d_\ell := d$ and $d_i, d_i^{\text{hid}} \in \N$ for $i=1,\ldots,\ell$. For the functions $\sigma_i$, we assume smooth (coordinate-wise) nonlinear activation functions, such as $\tanh$ or a sigmoid function, applied to each coordinate. Our vector field assumptions emphasize applications and are inspired by typical deep neural networks. As such, the assumptions to obtain unique solutions to \eqref{eq:node} are fulfilled but could be greatly relaxed. It is sufficient to assume $L^2$ integrable parameter functions and general vector fields with adequate smoothness in both variables \cite{LongoNovoObaya}. Minimal requirements to guarantee the existence of both the LEs and its gradients used in Section~\ref{sec:ftlereg} are an interesting question postponed to future work.
 
 The simplest examples of vector fields of interest in \eqref{eq:node} are
 \begin{equation*}
 	f(\theta, x) = \sigma(Wx+b),
 \end{equation*}
  with $\theta = [W,b]$ containing a weight matrix $W\in \R^{d\times d}$ and a bias vector $b\in \R^d$. Following the definitions from above, this vector field has only one parameter layer, $K=1$, and one vector field layer $\ell = 1$.
  
  \begin{remark}\label{rem:resnetnode}
  	The connection of the parameter layers defined here of nODEs to the layers of discrete Residual Neural Networks (ResNets) \cite{resnet2016} is as follows: Choose $K=T$ with equidistant $\beta_k-\alpha_k = 1$ and discretize the nODE \eqref{eq:node} with an explicit Euler scheme of step size 1. Then, for an input $x(0) = x_0$ and $x_k:=x(k)$ we obtain the typical ResNet update rule
  	\begin{equation}\label{eq:resnet}
  		x_{k} = x_{k-1} + f(\theta(\alpha_k), x_{k-1}) = x_{k-1} + f(\theta_{k}, x_{k-1}),\quad k=1,\ldots, K,
  	\end{equation}
  	which is a ResNet with $K$ layers \cite{weinan2017}. Let us also point out that in the general context of nODE implementations, it is not clear how to define the concept of ``nODE layers'' for solvers which produce variable discretization steps per forward pass \cite{CRBD18}.
  \end{remark}
 
 Neural ODEs of the form \eqref{eq:generalf} are fairly general and can realize widely different dynamics and, as a result, very different discrete network architectures. To see this, let us assume the total number of real-valued scalar parameters included in $\theta(t)$ of form \eqref{eq:theta} to be fixed.
 \begin{itemize}
 	\item  If we fix the number of vector field layers to one, $\ell =1$, but a large number of parameter layers $K$, the dynamics correspond to deep but narrow neural network architectures. The potential complexity of the input-output flow is realized via the time variable. Rather elementary piecewise-autonomous dynamics are combined in time sequence to reach high expressivity.
 	
 	\item If there is just one parameter layer, $K = 1$, but the number of vector field layers $\ell$ is large, the input-output structure is similar to Recurrent Neural Networks (RNNs) \cite{pascanu2013, engelken2020lyapunov}. The potential complexity of the input-output flow is realized instantaneously in the phase-space variable. While the flow is autonomous, the model's expressivity is achieved through the vector field's layered structure.
 \end{itemize}
 
While highly relevant in applications, the two cases described above represent only the edge cases of \eqref{eq:generalf}. Combining different numbers of vector field layers and parameter layers leads to a large variety of possible models. For instance, in image recognition, ResNets of the form \eqref{eq:resnet} typically employ multi-layer transformations $f$, such as convolutional and pooling layers. Additionally, a further relevant modeling choice is given by the discretization procedure (necessary for any implementation) that approximates the continuous flow.

It is an important task to investigate which model architectures are advantageous for different objectives. The primary criteria are prediction accuracy, computational efficiency and trainability. However, while these are well understood, secondary criteria such as adversarial robustness, stability and interpretability still require improvements.

\subsection{Contribution}

We introduce \emph{finite-time Lyapunov exponents} (FTLEs) into the setting of neural ODEs \eqref{eq:node}. We demonstrate for classification problems that FTLEs are powerful indicators to investigate the behavior across very different deep neural network architectures. Supported by these insights, we present a novel training algorithm based on suppressing high FTLEs in order to reduce model complexity, input sensitivity and to improve adversarial robustness. For a qualitative understanding, we implement and visualize these methods for a toy data set in two dimensions.

\subsection{Related works and structure}

The work \cite{storm2023finite} studies FTLEs in the discrete setting of feedforward neural networks (or multi-layer perceptrons) of various widths and depths which stay constant across all hidden layers. The key insights are: decision boundaries can be linked to high maximum FTLE ridges. If the networks reach a certain width, they exhibit positive FTLEs everywhere. The reason is that the final affine linear layer is sufficient for a highly accurate separation task. It is realized via a projection from the sufficiently high-dimensional augmented space down to the label dimension. As a result, the values of the network weights up to the final layer are insignificant to the task and can be frozen before training. 


In \cite{engelken2020lyapunov}, the authors use gradient flossing for improved trainability for RNNs. The idea is that during the training process of recurrent neural networks, the modulus of the gradients is reduced periodically via a Lyapunov exponent regularization term. This reduces the modulus of the parameter gradients and prevents vanishing or exploding gradients. The work \cite{vogt20} also uses LEs to investigate the stability of training dynamics of RNNs. More recently, \cite{vockmeisel25} has studied the LE spectrum across different DNN settings and architectures. It finds that good model performance is linked to dynamics close to criticality, which is signaled by LEs close to the critical value $0$. In addition, the work also provides different examples of modified training that manipulates the LE values.

In Section~\ref{sec:ftledef}, we define Lyapunov exponents in both continuous and discrete settings and link them to Lagrangian Coherent Structures. In Section~\ref{sec:experiments}, we track finite-time LEs of different nODEs which solve a two-dimensional classification task. In Section~\ref{sec:ftlereg}, we introduce a modified training algorithm to improve adversarial robustness by suppressing large FTLEs. In summary, this provides a very detailed, yet transparent, perspective on the use and utility of FTLEs for nODEs on the level of phase space dynamics as well as regarding their training.

\section{Finite-time Lyapunov exponents}\label{sec:ftledef}

In this section, we briefly introduce Lyapunov exponents, both the theoretical aspects as well as their numerical implementation relevant for our deep learning setting. We then connect these quantities to Lagrangian coherent structures (LCS), linked to local maxima of the FTLEs.


\subsection{Lyapunov exponents}
Let us consider a general ODE initial value problem of the type \eqref{eq:node} with a vector field that for all $\theta$ satisfies $f(\theta, \cdot) \in C^1(\R^d, \R^d)$ and solution $x(t)$. Let the corresponding time-dependent flow be denoted as
\begin{equation}\label{eq:flow}
	\Phi(t,x_0) := x(t).
\end{equation}
Further, let the spatial Jacobian of the vector field $f(t,x)$ evaluated at position $x$ at time $t\geq 0$ be denoted as $\textnormal{D}_x f(t,x)$, i.e., it is the $d\times d$ matrix of partial derivatives of the vector field. Then the solution $Y(t,x_0)$ to the linearized vector field along a trajectory $(t,x(t))$ gives us the so-called \emph{variational equation}
\begin{equation}\label{eq:tangentode}
\begin{cases}
	\dot{Y}(t) = \textnormal{D}_x f(t,x(t)) Y(t), & t\in [0,T),\\
 Y(0)= \textnormal{Id} \in \R^{d\times d},
\end{cases}
\end{equation}
where $\textnormal{Id}$ denotes the identity matrix. The variational equation describes the tangent mapping of the flow along a solution $x(t)$. Note carefully that the solution of the variational equation implicitly depends also on the initial value $x_0$ of the trajectory $x(t)$ of the nODE. Therefore, we shall also sometimes write $Y(t,x_0)$ to emphasize this dependence. 
For a given initial value $x_0\in \R^d$ to \eqref{eq:node}, one can formally denote the tangent mapping by
\begin{equation}\label{eq:tangentflow}
	\textnormal{D}_{x_0} \Phi(t,x_0):= Y(t,x_0): \textnormal{T}_{x_0} \mathcal{M} \to \textnormal{T}_{\Phi(t,x_0)}\mathcal{M},
\end{equation}
where $\textnormal{T}_{x_0}\mathcal{M}$ denotes the tangent space to a manifold $\mathcal{M}$ at a point $x_0\in\mathcal{M}$. Since we only work in $\R^d$, we can always think of $\mathcal{M}$ as a submanifold of $\R^d$. We see from \eqref{eq:tangentode} that computing $Y(t,x_0)$ at some time point $t=t_1$ requires the knowledge of the whole trajectory $x(t)$ for $t\in[0,t_1]$.
Let the singular values of $Y(t,x_0)\in \R^{d\times d}$ be denoted in monotonically decreasing order as 
\begin{equation*}
\Lambda_1(t,x_0)\geq \ldots \geq \Lambda_d(t,x_0)> 0,
\end{equation*}
which are all positive as $Y(t,x_0)$ has full rank for $C^1$ vector fields $f$. Then, the $j$-\emph{th finite-time Lyapunov exponent (FTLE) of} $\Phi(t,x_0)$ is defined as
\begin{equation*}
	\lambda_j(t,x_0) :=  \frac{1}{t}\ln[\Lambda_j(t,x_0)],\quad t>0,\quad j = 1,\ldots, d,
\end{equation*}
and expresses the \emph{time-averaged logarithmic growth rate} of $Y(t,x_0)$ on $[0,t]$ along the $j$-th principal axis (defined via the singular value decomposition). The \emph{maximum FTLE}, 
\begin{equation}\label{eq:deflmax}
\lambda_{\max}(t,x_0):= \lambda_1(t,x_0) = \frac1t \ln \|\textnormal{D}_{x_0} \Phi(t,x_0)\|_2,
\end{equation}
measures the (linearized) maximum exponential divergence of nearby trajectories over the time interval $[0,t]$, where $\|\cdot\|_2$ denotes the Euclidean norm. 

\begin{remark}
Because of the convergence properties of the Lyapunov exponents, when considering infinite-time Lyapunov exponents, the initial conditions for the variational equation are of little importance \cite{oseledec68}, and $Y(0) = \textnormal{Id}$, with $\textnormal{Id}$ the identity matrix is widely accepted. For the finite-time Lyapunov exponents on the other hand, different initial conditions result in (potentially significant) different results, in particular in cases in which the integration time is short with respect to the intrinsic expansions and contractions of the system. Yet, it is standard to use the identity matrix also as a practical first benchmark case for FTLEs.
\end{remark}

\begin{remark}
For notational simplicity, we have so far restricted our discussion to initial conditions at $t = 0$. More generally, when computing FTLEs over an arbitrary time interval $[t_0, t_1]$ with $x(t_0) = x_0$, we write $\lambda_j([t_0, t_1], x_0)$ for the $j$-th FTLE.
\end{remark}

\subsection{Discretized setting}
While our general setting is a continuous nODE \eqref{eq:node}, any numerical implementation requires a discrete solver of the dynamics and hence is realized in a discrete setting. As an important example (cf.\ Section~\ref{sec:experiments}), we write down the Euler discretization  with step size $\Delta t>0$ of \eqref{eq:node} (cf.\ Remark~\ref{rem:resnetnode} which used step size $\Delta t = 1$). This yields the update rule
\begin{equation*}
\Phi(t_n,x_0) := 	x_{n} = x_{n-1} + \Delta tf(t_{n-1},x_{n-1}),\quad n=1,\ldots ,N,
\end{equation*}
denoting $t_n:=n\Delta t$, $x_n:=x(t_n)$. To obtain the tangent mapping \eqref{eq:tangentflow} in the discrete setting, one applies the iterative chain-rule
\begin{align*}
\textnormal{D}_{x_0} \Phi(t_n,x_0) &= \frac{\partial x_{n}}{\partial x_{n-1}}\frac{\partial x_{n-1}}{\partial x_{n-2}}\cdots \frac{\partial x_{1}}{\partial x_{0}}\\
& = (\textnormal{Id} + \Delta t \textnormal{D}_{x} f(t_n,x_n))\cdots (\textnormal{Id} + \Delta t \textnormal{D}_{x} f(t_0,x_0)),\numberthis \label{eq:chainrule}
\end{align*}
where $\frac{\partial x_n}{\partial x_{n-1}}$ denotes the one time-step Jacobian. Then, denoting the singular values of  $\textnormal{D}_{x_0}\Phi(t_n,x_0)$ as  $\Lambda^n_1(x_0)\geq \cdots \geq \Lambda^n_d(x_0) > 0$, the discrete FTLEs at mapping $n=1,\ldots, N$ are defined as
\begin{equation}
		\lambda^n_j(x_0) :=  \frac{1}{t_n}\ln[\Lambda^n_j(x_0)],\quad n=1,\ldots,N; j = 1,\ldots, d.
\end{equation}

Our setting is not restricted to an Euler discretization of \eqref{eq:node} and  allows for any discretization scheme. Then, \eqref{eq:chainrule} needs to be adapted accordingly. For instance adaptive solvers for the numerical implementation can vary their step size in dependence on a predefined error tolerance and model application (training, inference, etc.) \cite{CRBD18}, see Section~\ref{subsec:ftleregeff}.

\subsection{Lagrangian coherent structures}\label{subsec:lcs}

The basic idea of coherent structures in dynamical systems is to identify subsets $\mathcal{S}\subset \R^d$ of the phase space, which stay ``coherent'' under the flow. Mathematically, the most classical definition works via invariant sets and is focused on autonomous systems. A set $\mathcal{S}$ is invariant under a flow $\phi$ if $\phi(t,\mathcal{S})=\mathcal{S}$ for all $t\in\R$~\cite{GH}. Yet, this definition is too rigid for many applications as it requires \emph{complete invariance} for \emph{all times}, and it requires a \emph{flow} $\phi$. To avoid full invariance, the most elegant modification is to consider almost invariant sets. For a fixed, small $\varepsilon>0$ and a final time $T>0$ for $t\in [0,T]$, we say that a set $\mathcal{S}$ is $\varepsilon$-coherent (or almost invariant) over time $T$ if 
\begin{align}
    \frac{\mu\left(\mathcal{S}\cap \phi(-T,\mathcal{S})\right)}{\mu(\mathcal{S})}\geq 1-\varepsilon ,
\end{align}
where $\mu$ is a reference measure taken usually as the Lebesgue measure. In particular, almost invariant sets are those where just a small amount of mass leaks under the flow~\cite{DellnitzJunge1,FroylandPadberg}. Even this definition is conceptually not yet sufficient as we still require a classical flow $\phi$, which is suitable for autonomous ODEs but our nODEs are typically non-autonomous. A natural alternative is to use a more general definition and consider the solution mapping 
\begin{equation*}
x(t_1)=\Phi(t_1,t_0;x_0)
\end{equation*}
of our neural ODE, which is solved over a time interval $[t_0,t_1]$, $t_0 < t_1$, with initial condition $x(t_0)=x_0$. To account for the non-autonomous setting, the almost invariant sets must be made time-dependent. Instead of just $\mathcal{S}$, we consider a time-dependent family of sets $\mathcal{S}_t$ for $t\in[t_0,t_1]$. One then defines the set in the beginning to be coherent with the set at the end of the time interval by considering conditions such as
\begin{align}\label{eq:tinvariant}
    \frac{\mu\left(\mathcal{S}_{t_0}\cap (\Phi^{-1}(t_1,t_0;\mathcal{S}_{t_1}))\right)}{\mu(\mathcal{S}_{t_0})}\geq 1-\varepsilon. 
\end{align}
Yet, computing every almost invariant set can be computationally challenging from these definitions. Similar to the situation in chaotic dynamical systems, one looks for simpler observables to detect structures that could be coherent. It turns out that FTLEs are linked to almost invariance in many (yet not all) situations, i.e., FTLE ridges can function as useful observables~ \cite{haller00,haller02,ShaddenLekienMarsden,haller11,NolanSerraRoss}. Here we adopt this viewpoint and work with FTLEs. To define them, consider the right Cauchy-Green stress tensor
\begin{align*}
 \mathcal{C}(t_1,t_0;x_0):=(\textnormal{D}_{x_0}\Phi(t_1,t_0;x_0))^\top \textnormal{D}_{x_0}\Phi(t_1,t_0;x_0).  
\end{align*}
Generically, $\mathcal{C}(t_1,t_0;x_0)$ is symmetric and positive definite with eigenvalues 
\begin{align*}
\rho_1 \geq \rho_2 \geq \cdots \geq \rho_d > 0
\end{align*}
with associated orthonormal eigenvectors $\xi_1,\ldots,\xi_d\in\R^d$. A repelling Lagrangian coherent structure (LCS) can then be defined as codimension one manifold $\mathcal{M}$ that locally maximize normal repulsion. More formally, this means that if $\vec{n}$ denotes a unit normal vector at $x_0$ to $\mathcal{M}$ one can define the stretching factor
\begin{align*}
 \rho_*(x_0,\vec{n}):=\sqrt{\vec{n}^\top \mathcal{C}(t_1,t_0;x_0) \vec{n}}   
\end{align*}
and require maximal normal stretching
\begin{align*}
 \rho_*(x_0,\vec{n})=\rho_1,\qquad \vec{n}=\xi_1.   
\end{align*}
In addition, it is common to require a spectral gap via the condition $\rho_{2}<\rho_1$ and a natural tangency condition
\begin{align*}
    \textnormal{T}_{x_0} \mathcal{M}=\textnormal{span}(\xi_d,\ldots,\xi_{2}).
\end{align*}
It is then standard to observe that the maximum Lyapunov exponents can also be expressed via the eigenvalues of the right Cauchy-Green stress tensor
\begin{align}\label{eq:cauchygreen}
    \lambda_{\max}([t_0,t_1],x_0)=\frac{1}{t_1 - t_0} \ln \rho_1,
\end{align}
where $\rho_1$ depends implicitly on $x_0$ and the time interval as is evident from its definition above. In particular, maximizing the top FTLE is equivalent to maximizing $\rho_1$, which locally is part of the definition of the LCS. Then, we define a ridge $\mathcal{R}$ of the FTLEs as a codimension-one manifold, where the top Lyapunov exponent is maximal in the direction normal to the ridge. In general, one expects an LCS $\mathcal{M}$ to be contained in the set of ridges, i.e., $\mathcal{M}\subset \mathcal{R}$. However, there are cases where equality might fail and regularity issues arise~\cite{KarraschHaller}. In some sense, this is not surprising as spectral observables such as FTLEs are always just a coarser observable than knowing all trajectories; this is similar to classical chaos detection, where a positive LE or FTLE is an indicator but not a proof of chaos. However, ridges of FTLEs are easier to compute and can often be efficiently estimated by standard algorithms making them quite practical tools, e.g., in fluid dynamics. As we shall demonstrate below, the same conclusion applies to nODEs and their training. 



\section{Tracking FTLEs of neural ODEs}\label{sec:experiments}
This and the following section address the numerical implementations of neural ODEs in a two-dimensional setting. We provide a detailed discussion for two representative minimal examples.

\subsection{Numerical setup}\label{subsec:numerics}
 Before defining the two nODE examples, we describe the numerical setup in which they are implemented\footnote{The implementation code is available at \url{https://github.com/twoehrer/FTLEs_of_nODEs.git}}.\\

\begin{figure}
		\centering 
		\begin{subfigure}{0.27\textwidth}
			\includegraphics[width=\textwidth, trim={0 0.em 0 0},clip]{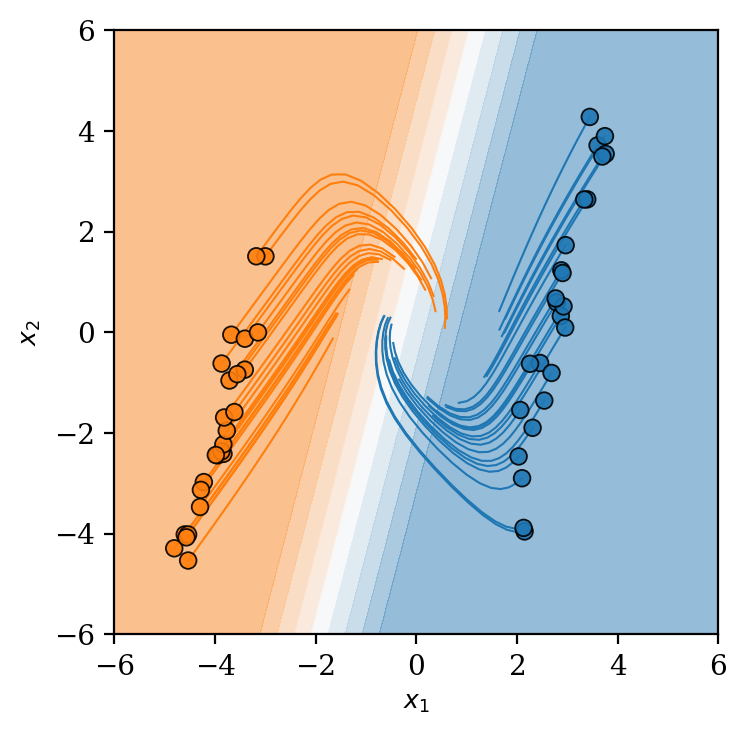}
			\caption{nODE trajectories \label{subfig:moons:traj1param}}
		\end{subfigure}
		\begin{subfigure}{0.35\textwidth}
			\includegraphics[width = \textwidth, trim={0 2.em 0 0},clip]{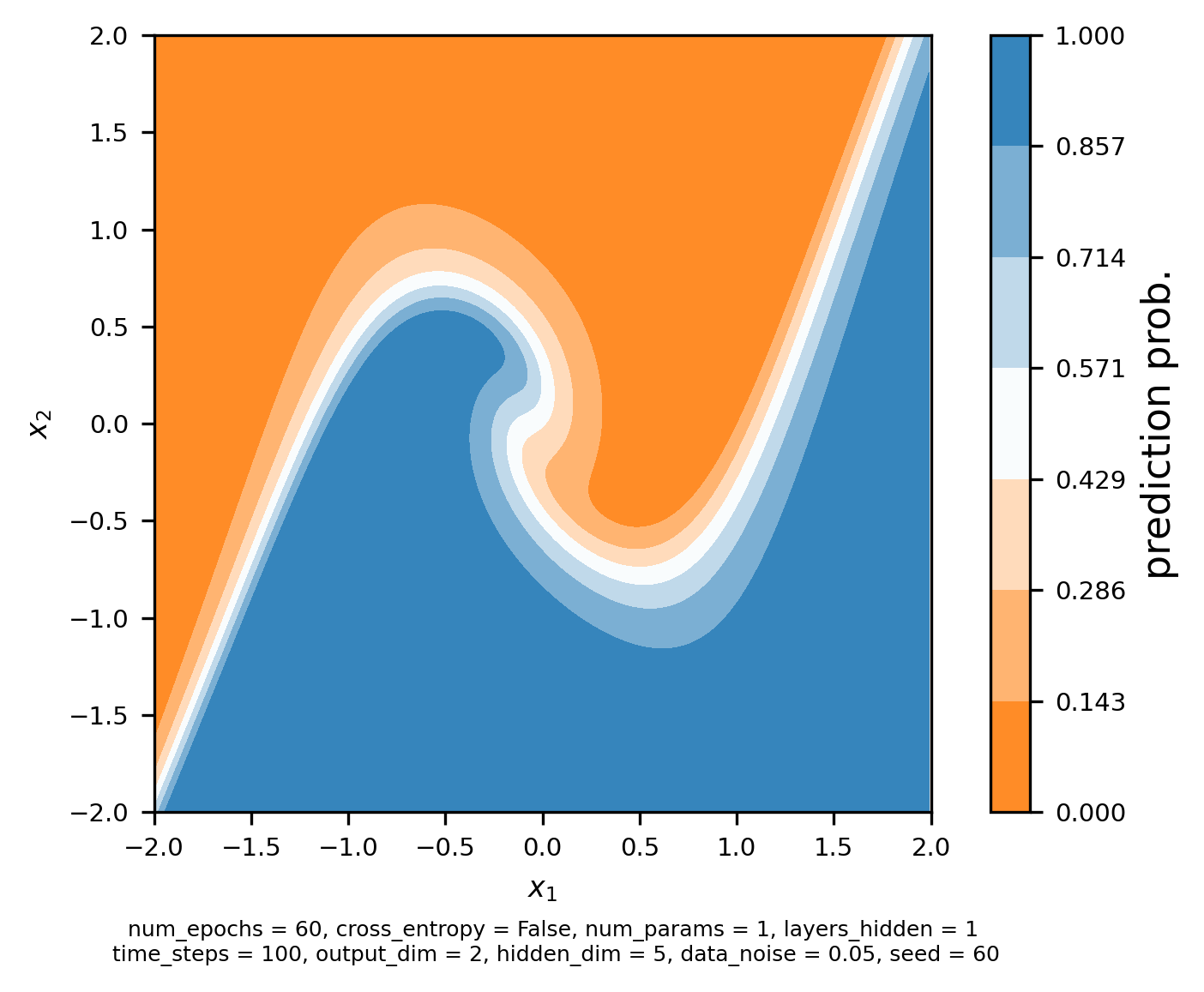}
			\caption{Prediction level sets \label{subfig:moons:levelsets1param}}
		\end{subfigure}\\
        \begin{subfigure}{0.3\textwidth}
			\includegraphics[width = \textwidth, trim={0 0.em 0 0},clip]{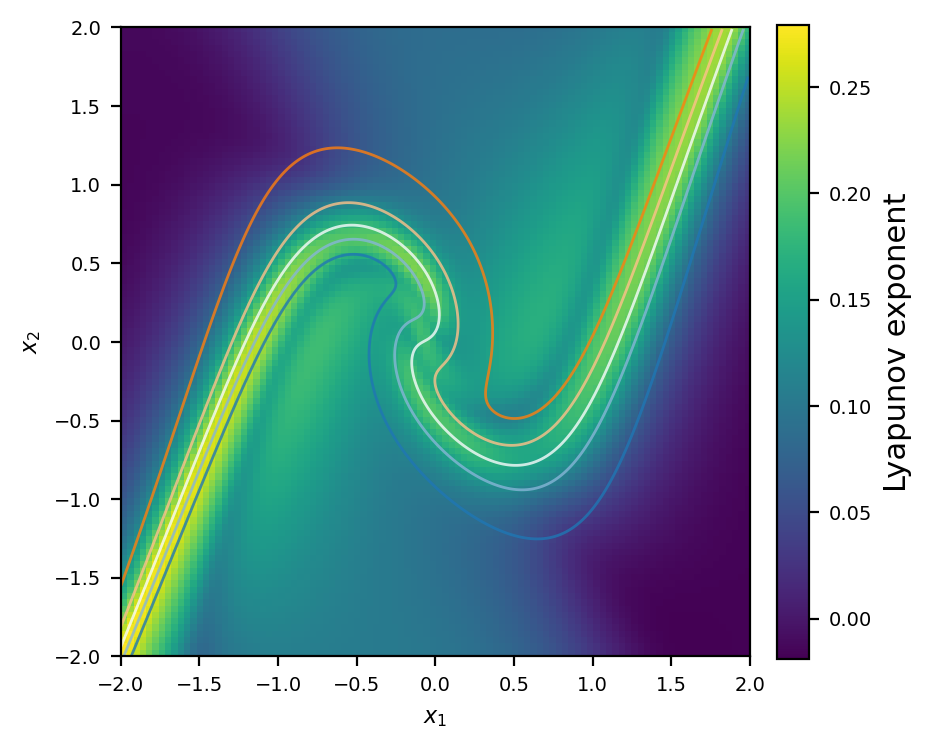}
			\caption{\centering max.\ FTLEs on [0,10]\\and level sets \label{subfig:moons:ftlesandlevel}}
		\end{subfigure}
		\begin{subfigure}{0.315\textwidth}
			\includegraphics[width = \textwidth, trim={0 0.em 0 0},clip]{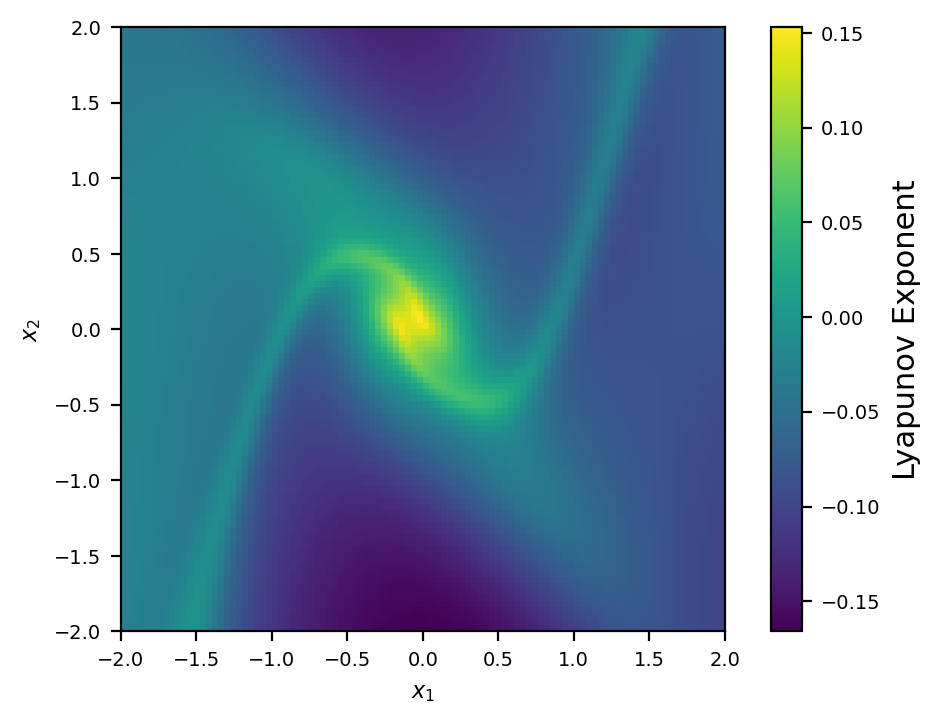}
			\caption{min.\ FTLEs on $[0,10]$ \\~\label{subfig:moons:minftles1param}}
		\end{subfigure}
		
		\caption{Plots for Example~\ref{ex:1param}}
	\end{figure}

\noindent\textbf{Discretization:}
Our numerical experiments specify the right-hand side of \eqref{eq:node}--\eqref{eq:generalf} with activation functions $\sigma_i(y) = \tanh(y)$. We utilize an explicit Euler discretization with step size $\Delta t = 0.1$ and final time $T=10$. The simple Euler method allows a straight-forward computation of the FTLEs. Reducing the approximation error further in test cases, via reduced step size or higher-order schemes (using the \texttt{torchdiffeq} package \cite{CRBD18}) has not led to any qualitative changes in the dynamics discussed in this work. The choice $T \gg 1$ allows for sufficient trajectory evolution relative to the boundedness of the vector field.\\
\textbf{Parameter training}:
The training dataset consists of $4000$ points of the moons dataset with colored labels (see Figure~\ref{subfig:trainingset}). We selected this specific dataset to avoid topological issues present in, for example, the circle toy datasets (cf.\ \cite{dupont2019}). We encode the class labels as vectors $y_{\text{blue}} = (0,1)$ and $y_{\text{orange}}= (0,-1)$ and minimize the \emph{empirical mean square loss} over a randomly selected batch of inputs and labels, $({\mathbf {x_0, y}}): = (x_{0,i}, y_i)_{i = 1}^{\Ndata}$, as
\begin{equation}\label{eq:optimization}
    \min_{\theta} \mathcal{L}({\mathbf {x_0, y}}, \theta)
\end{equation}
with loss function
\begin{equation}\label{eq:loss}
    \mathcal{L}({\mathbf {x_0,y, \theta}}):=\frac{1}{\Ndata}\sum_{i=1}^{\Ndata} (L\Phi(T,x_{0,i})-y_i)^2.
\end{equation}
The output layer function $L:\R^2 \to \R^{2}$ is an affine linear layer with trainable parameters and $\Phi(T,\cdot)$ denotes the time-$T$ flow generated by the neural ODE. We opted for two-dimensional outputs and labels (compared to prescribing a 1d projection of $L$ and 1d labels) in order to minimize the separation role of the output projection $L$ in relation to the nODE flow (cf.\ \cite{storm2023finite}). The nODE parameters are initialized via standard He initialization \cite{resnet2016}. The parameter optimization \eqref{eq:optimization} is realized via the stochastic gradient descent based Adam algorithm \cite{kingma14adam} and $\Nbatch = 64$. The parameter gradients are computed via Pytorch's standard auto-differentiation capabilities. The initialization is consistent across all plots. A further modification of the loss, based on an added FTLE suppression term, is introduced in Section~\ref{sec:ftlereg}.\\
\textbf{Model predictions and decision boundary:} The model classification depends on two components. First, the affine linear map $L$ introduced in \eqref{eq:loss} defines a partition of the state space. Second, the nODE flow $\Phi_T$ maps inputs into these specific regions. The model's classification probabilities are then based on the relative Euclidean distances of $L\Phi_T(x_0)$ to the class labels $y_\text{blue}$ and $y_\text{orange}$, see Figure~\ref{subfig:moons:traj1param}. We denote $\operatorname{pred}(x_0) \in [0,1]$ as the resulting model prediction confidence of an input $x_0$ to be classified with label $y_\text{blue}$ (and hence $1-\operatorname{pred}(x_0)$ describes the confidence for $y_\text{orange}$). This means if a model assigns an input $x_0$ a prediction that satisfies $\operatorname{pred}(x_0) > 0.5$, it is classified as $\text{blue}$ and if $\operatorname{pred}(x_0) < 0.5$ as $\text{orange}$. As a consequence, the \emph{decision boundary} (or $\eps$-decision margin) is the level set of $\eps$-low confidence in both labels,
\begin{equation}\label{eq:Deps}
    D_\eps:=\{x_0\in\R^2 : |\operatorname{pred}(x_0) - 0.5|< \eps \}, \quad \eps\geq 0.
\end{equation}
The choice $\eps = 0$, describes the level set of inconclusive predictions which describes the boundary of the class subsets. In practice (e.g., the white strip in Figure~\ref{subfig:moons:levelsets1param}), one considers an $\eps>0$ tolerance that determines when an input is deemed undecided by the model.\\

    \subsection{Vector field examples in two dimensions}
\begin{figure}
		\centering 
		\begin{subfigure}{0.32\textwidth}
			\includegraphics[width=\textwidth, trim={0 0.5em 0 0},clip]{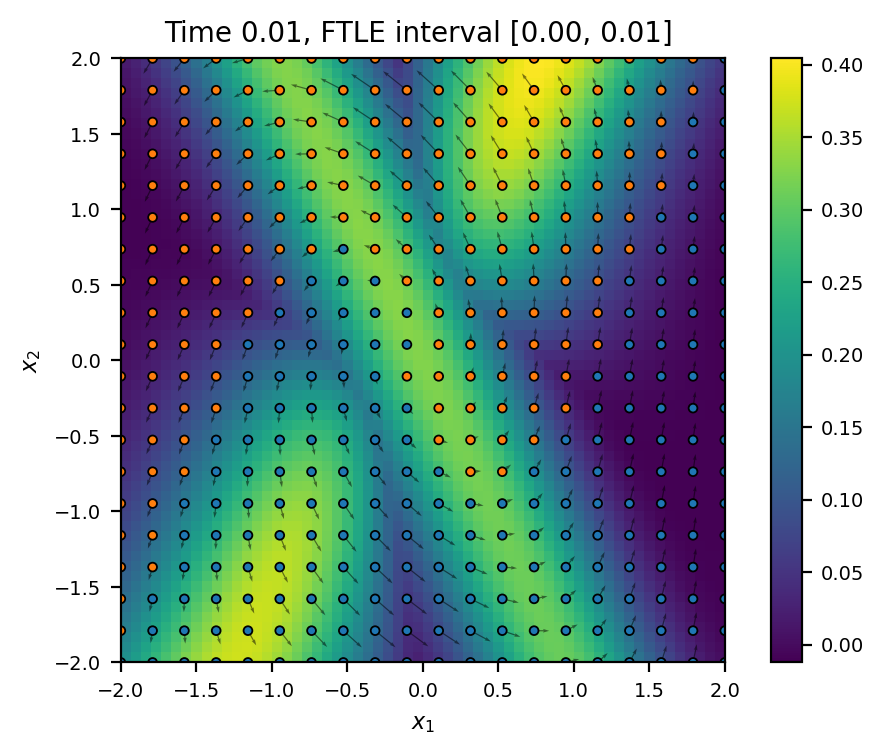}
		\end{subfigure}
		\begin{subfigure}{0.32\textwidth}
			\includegraphics[width=\textwidth, trim={0 0.5em 0 0},clip]{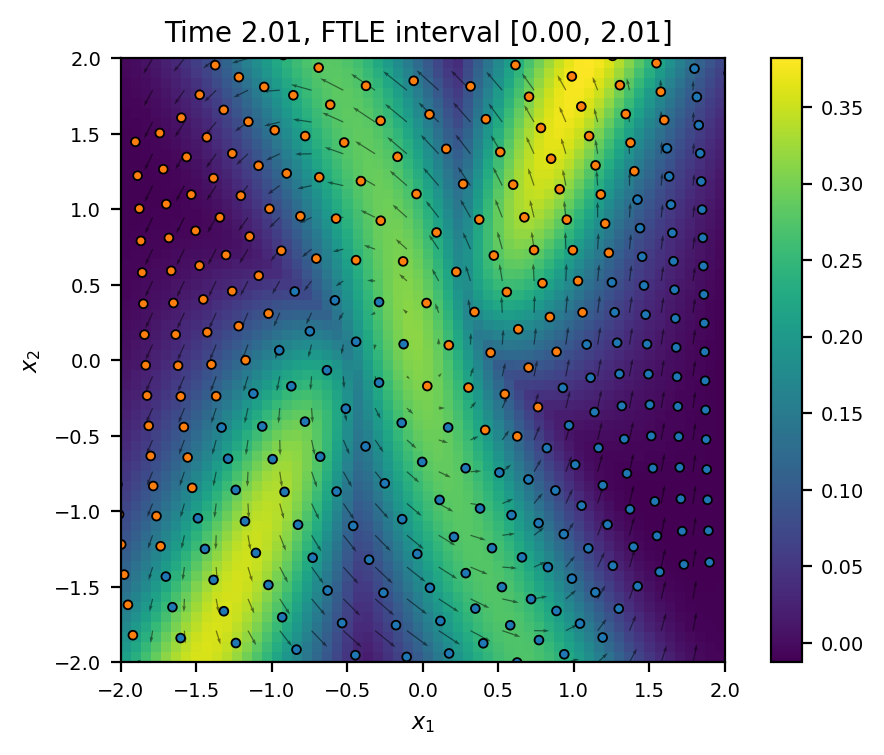}
		\end{subfigure}
		\begin{subfigure}{0.32\textwidth}
			\includegraphics[width=\textwidth, trim={0 0.5em 0 0},clip]{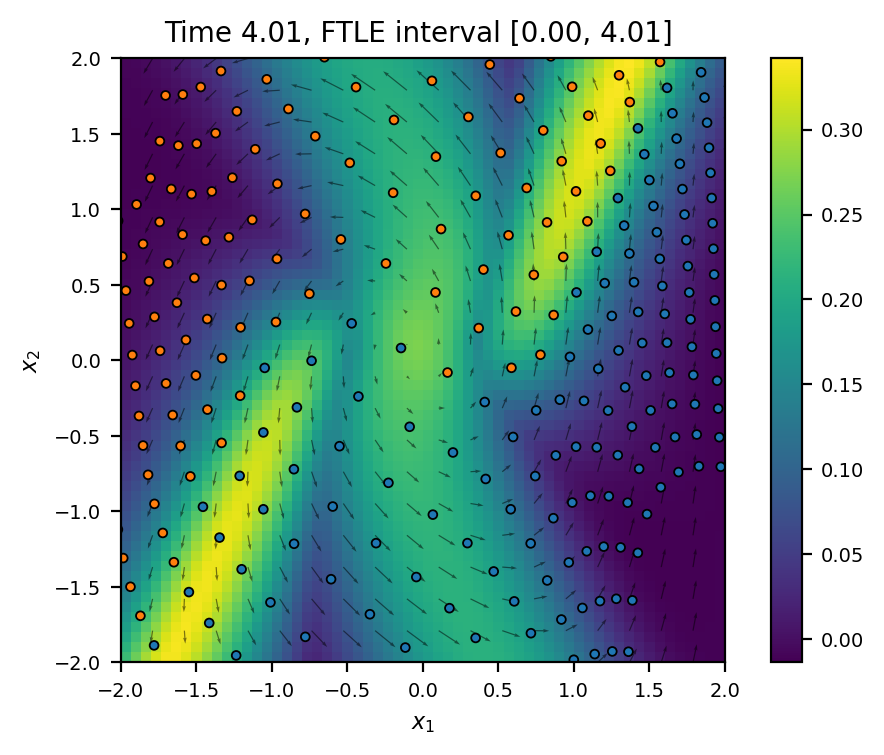}
		\end{subfigure}
		\\
		\begin{subfigure}{0.32\textwidth}
			\includegraphics[width=\textwidth, trim={0 0.5em 0 0},clip]{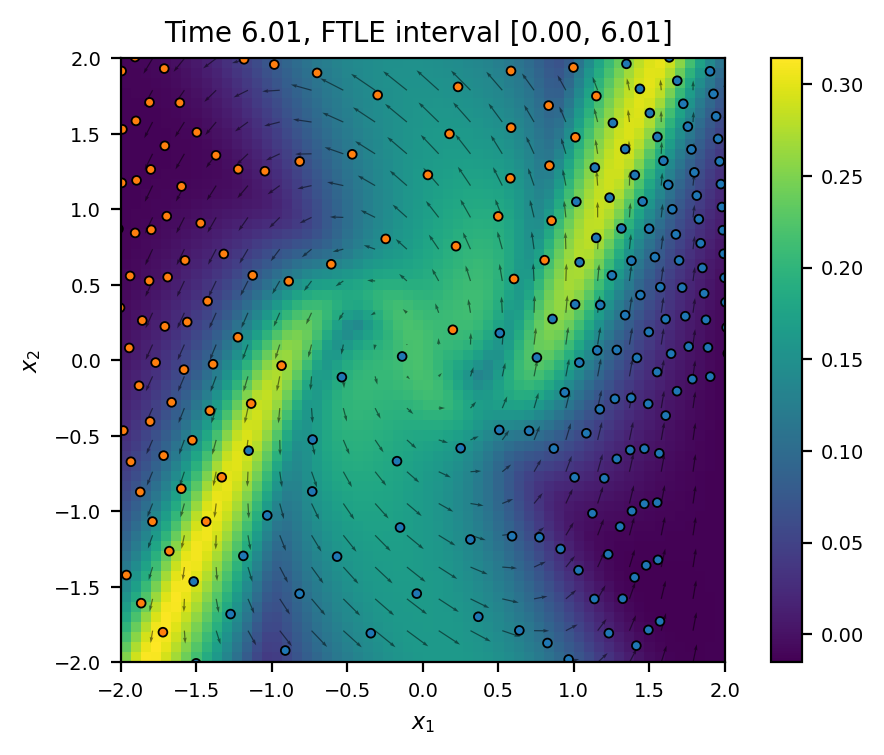}
		\end{subfigure}
		\begin{subfigure}{0.32\textwidth}
			\includegraphics[width=\textwidth, trim={0 0.5em 0 0},clip]{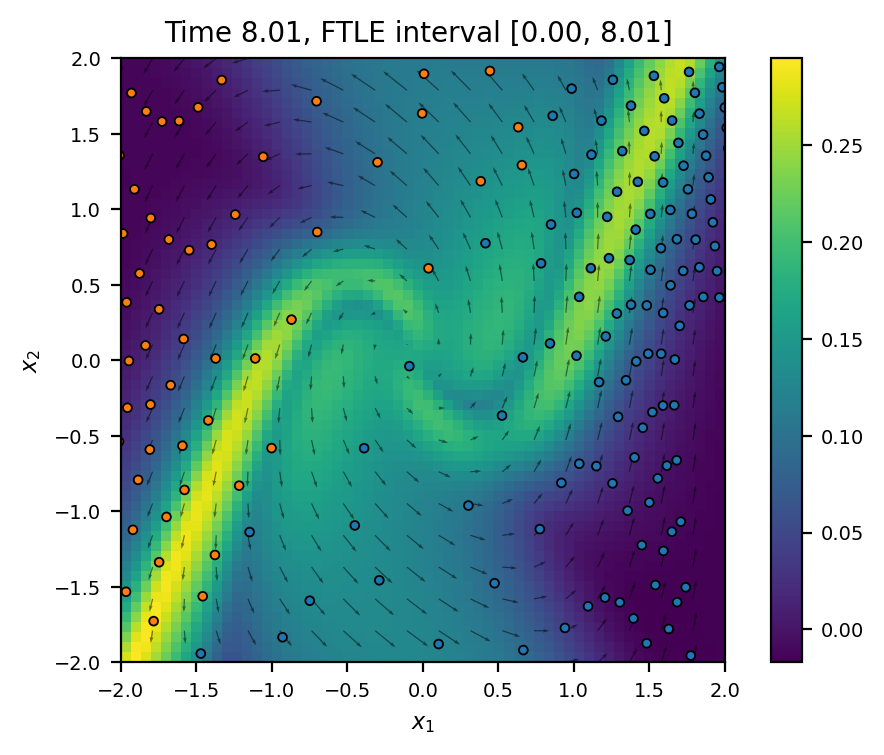}
		\end{subfigure}
		\begin{subfigure}{0.32\textwidth}
			\includegraphics[width = \textwidth, trim={0 0.5em 0 0},clip]{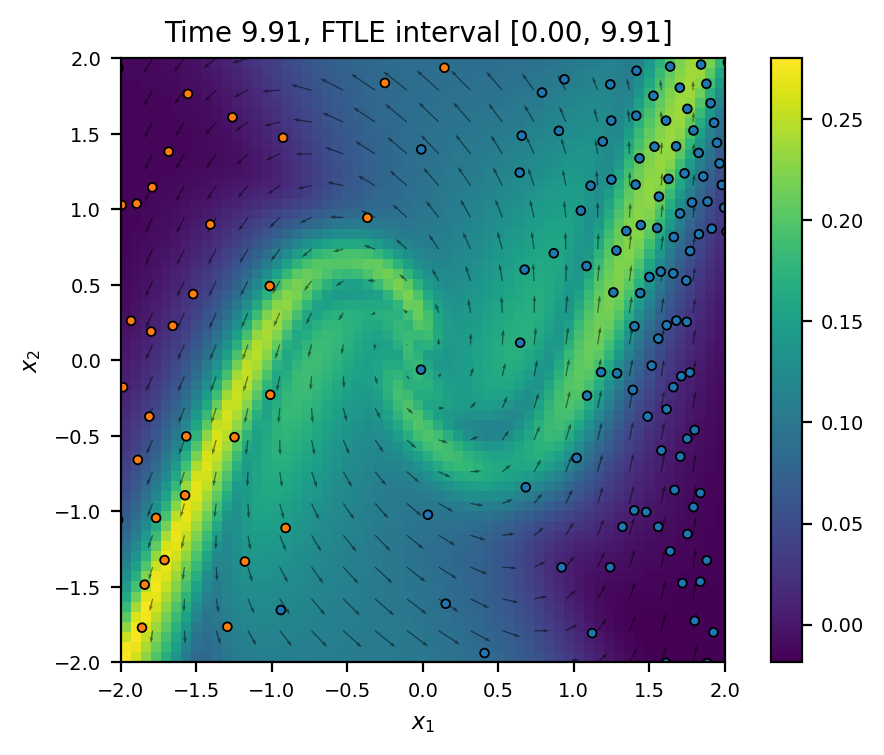}
		\end{subfigure}
		
		\caption{\small Plotted are the neural ODE evolution of Example~\ref{ex:1param} for a grid of initial points with labels colored orange/blue according to the model's final prediction. The background colors depict the FTLEs for each point in $[-2,2]^2$ from initial time to current time. This results in an \emph{increasing FTLE interval} as time progresses. \label{fig:growingFTLE1param}}
	\end{figure}

	\setexample{1}
	\begin{subexample}
		\label{ex:1param}
			The first example of a neural ODE of type \eqref{eq:generalf} we consider has only one parameter layer, $K=1$ in \eqref{eq:theta}, resulting in an autonomous vector field $f(\theta, x)$. For sufficient expressivity, we use two vector field layers,  $\ell = 2$ in \eqref{eq:generalf} and the hidden layer width $d^{\text{hid}}_1 = d^{\text{hid}}_2 = 5$. Specifically,
		\begin{equation}\label{eq:vectorfieldex1}
			f(\theta, x) =  f_2 \circ f_1(x) = V_2\sigma(W_2\sigma (W_1x + b_1) + b_2) + a_2, \quad x\in\R^2,
		\end{equation}
		with constant-in-time parameters
		\begin{equation*}
			W_1 \in \R^{5\times 2}, b_1 \in \R^5, W_2 \in \R^{5 \times 5}, b_2 \in \R^5, V_1 = \textnormal{Id}\in \R^{5\times 5}, V_2 \in \R^{2\times 5}, a_2 \in \R^2.
		\end{equation*}
		
		As described in Section~\ref{subsec:numerics}, we consider the trained nODE parameters on the moons dataset of Figure~\ref{subfig:trainingset}. This leads to the model's prediction level sets and decision boundary of Figure~\ref{subfig:moons:levelsets1param}. For this fairly simple and well separated dataset the model achieves close to perfect classification accuracy on standard test data. Our main goal for this section is to analyze the model's dynamics that realize the classification result with the help of FTLEs.
		
Since the vector field \eqref{eq:vectorfieldex1} is autonomous and continuous, trajectories are unique and cannot intersect. This results in the qualitative behavior shown in Figure~\ref{subfig:moons:traj1param}, where the spatial flow is more constrained than in the non-autonomous case of Example~\ref{ex:5params} below. However, this does not prevent distinct inputs $x_0 \neq \tilde{x}_0$ on the same trajectory from being assigned different classes, as predictions depend only on the state at finite time $T$. For example, one might have $\Phi(t_1, x_0) = \tilde{x}_0$ for some $t_1 > 0$ while $\operatorname{pred}(x_0) \neq \operatorname{pred}(\tilde{x}_0)$.
 This fact highlights why finite-time Lyapunov exponents are our preferred perspective here over the $t \to \infty$ time-asymptotic dynamics. 

Our first observation regards the input trajectories of Figure~\ref{subfig:moons:traj1param}. The trained model achieves the separation by generating two qualitatively different dynamical phases for the training data. Near the origin, trajectories spiral outward anti-clockwise, effectively contracting points of the same class (both blue and orange) towards each other. Further away from the origin, the trajectories move approximately linearly and parallel, pushing inputs of different classes in opposite directions.
		
In Figures~\ref{subfig:moons:ftlesandlevel}--\ref{subfig:moons:minftles1param} (as well as Figure~\ref{subfig:intronoreg} from Section~\ref{sec:intro}), we plot the FTLEs of the neural ODE flow $x_0 \mapsto \Phi(T,x_0)$ for a grid of inputs $x_0 \in [-2,2]^2$. Our main observation is that the shape of the bright FTLE ridge of both the maximal FTLEs (Figure~\ref{subfig:moons:ftlesandlevel}) and the minimal FTLEs (Figure~\ref{subfig:moons:minftles1param}) captures the shape of the decision boundary. This signifies that the dynamics play an integral part in the separation task at hand (cf.\ \cite{storm2023finite}). \textbf{This shows that maximum FTLE ridges (computed on the whole time interval) are closely linked to the decision boundaries.} See Figure~\ref{subfig:moons:ftlesandlevel} for an overlay of the FTLE ridge, the decision boundary and further prediction level sets. One observes that the ridge and the decision boundary do not fully overlap. This is not due to numerical errors. They arise because the FTLE excludes the final affine-linear layer $L$, whereas $L$ is included in $\operatorname{pred}(x_0)$.

Figure~\ref{fig:growingFTLE1param} presents the FTLEs computed on growing subintervals of the full interval. It decodes how the decision boundary is generated dynamically in detail. While the underlying dynamics are time-independent, it takes about half the evolution time for the nonlinear structure of the decision boundary to appear as an FTLE ridge. The short-time plots exhibit the largest LE values and reveal that the input's early dynamics dominate the linearization process of the class separation. We refer to Section~\ref{subsec:ftlereg}, where we apply this knowledge to improve training. Since the process takes place across several stages in time, we conjecture that a detailed future analysis could rely on a multi-scale dynamics perspective: In a first fast phase, the dynamics generate manifolds that serve as a rough linear approximation of the class boundary before converging in a slow phase to a manifold which realizes the precise classification.
	\end{subexample}
	
	\newpage
	
\begin{figure}
	\centering 
	\begin{subfigure}{0.34\textwidth}
		\includegraphics[width = \textwidth, trim={0 2.8em 0 0},clip]{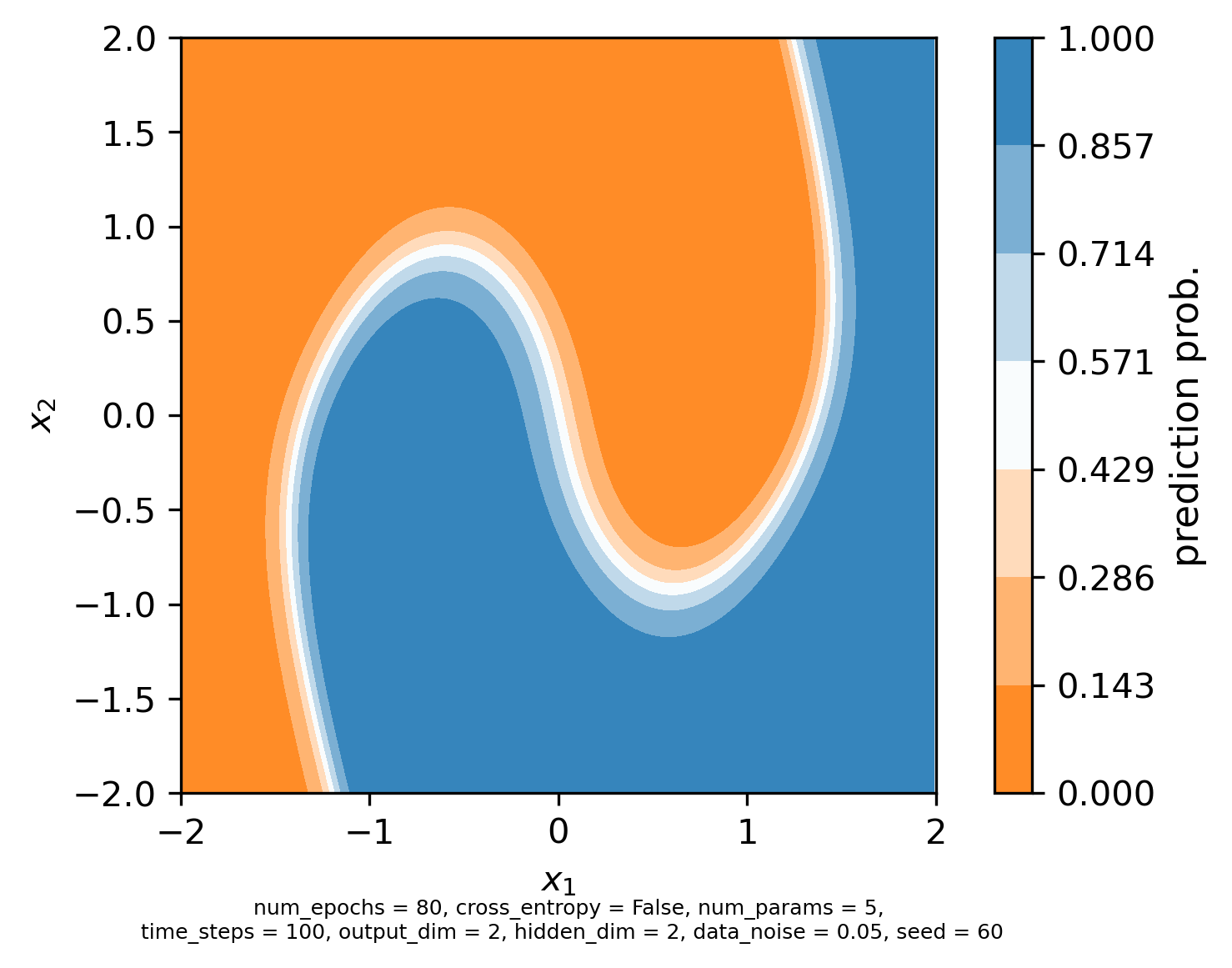}
		\caption{Prediction level sets \label{subfig:moons:levelsets}}
	\end{subfigure}
	\begin{subfigure}{0.32\textwidth}
		\includegraphics[width = \textwidth, trim={0 0.em 0 0},clip]{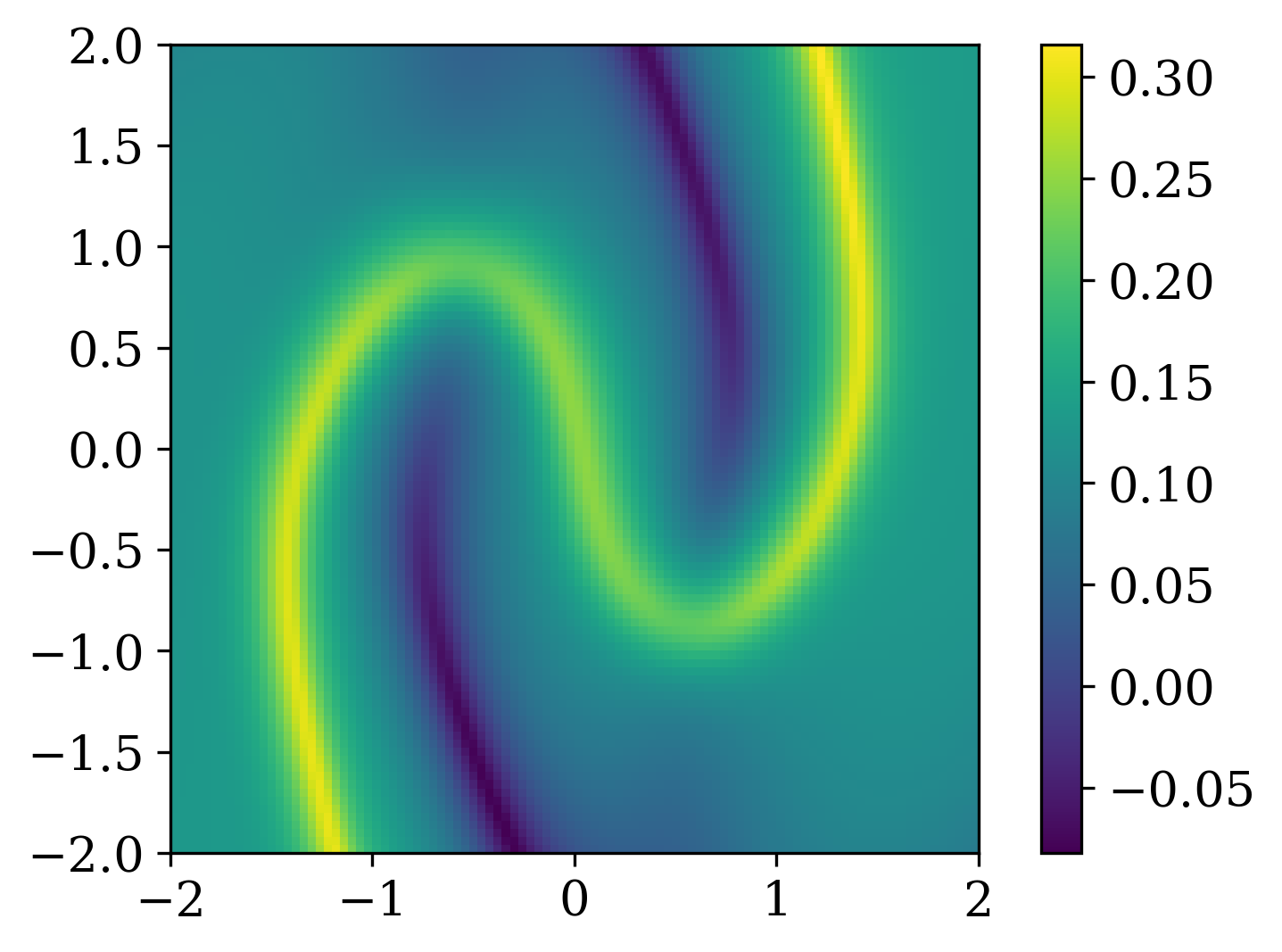}
		\caption{max.\ FTLEs on $[0,10]$ \label{subig:moons:ftles}}
	\end{subfigure}
	\begin{subfigure}{0.32\textwidth}
		\includegraphics[width = \textwidth, trim={0 0.em 0 0},clip]{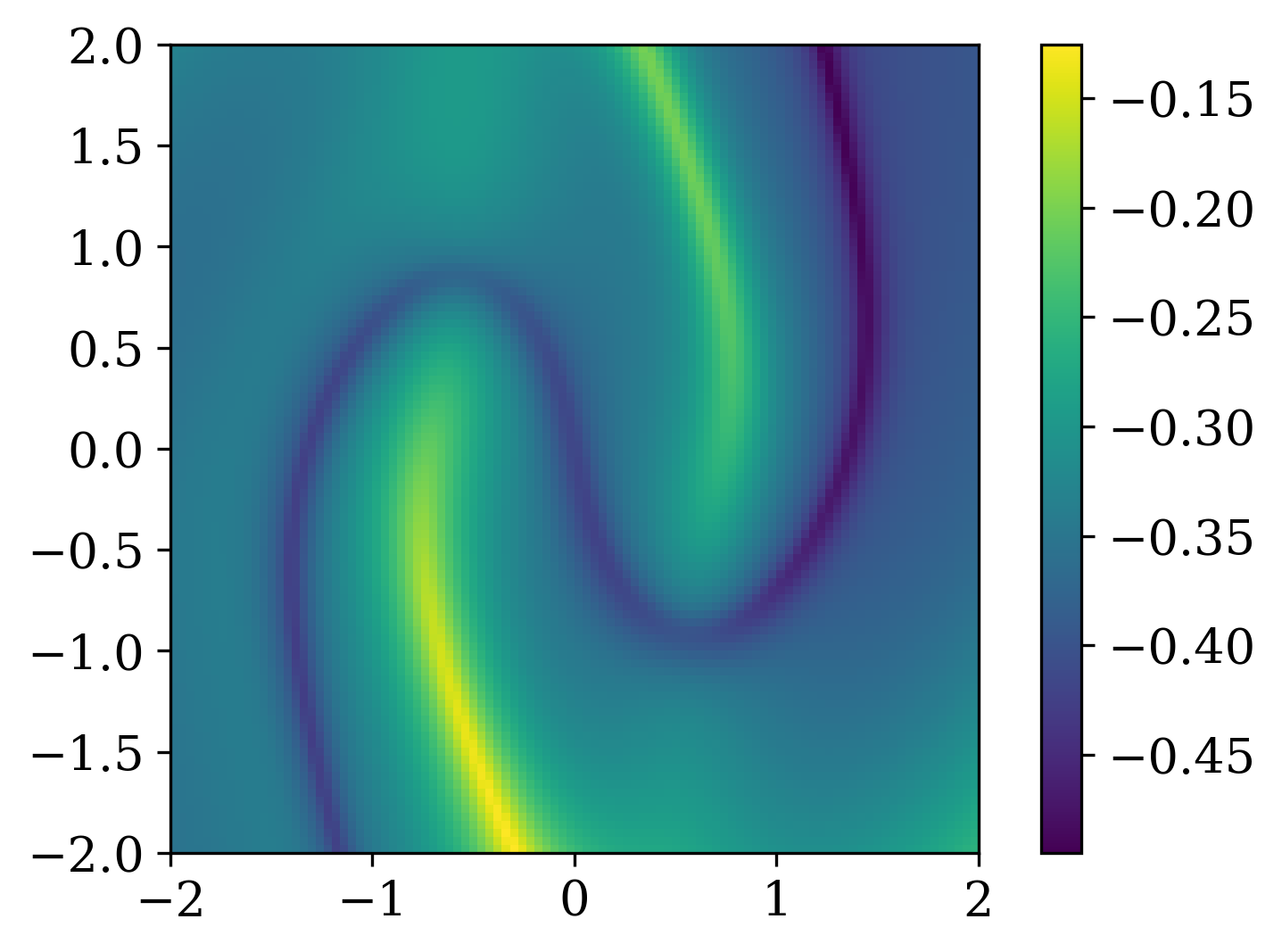}
		\caption{min.\ FTLEs on $[0,10]$ \label{subig:moons:minftles}}
	\end{subfigure}
	\caption{\small Plots for Example~\ref{ex:5params} \label{fig:5params}}
\end{figure}

\begin{figure}
	\centering 
	\begin{subfigure}{0.32\textwidth}
		\includegraphics[width=\textwidth, trim={0 0.5em 0 0},clip]{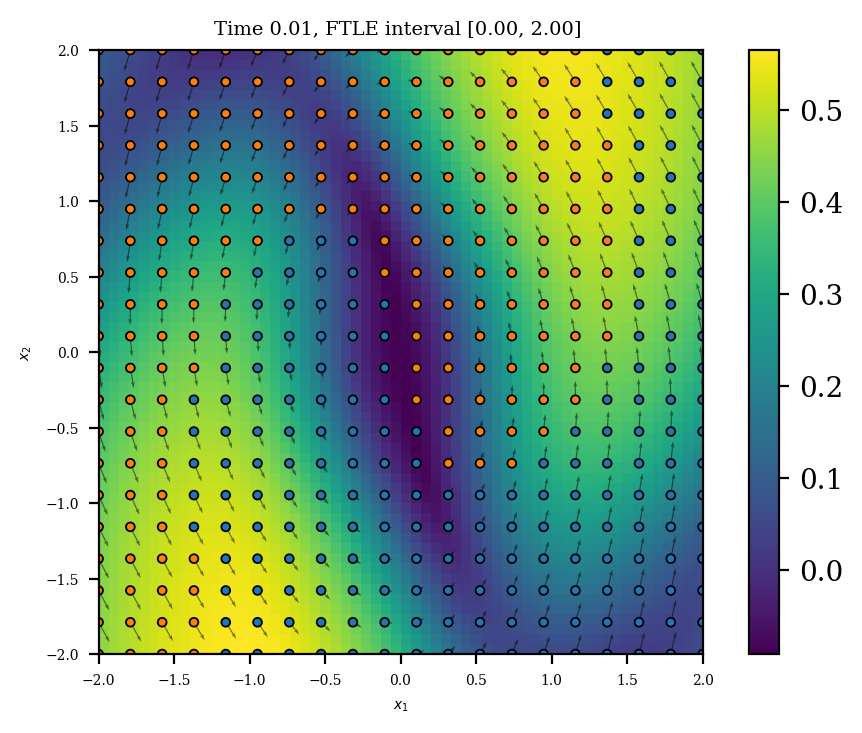}
	\end{subfigure}
	\begin{subfigure}{0.32\textwidth}
		\includegraphics[width = \textwidth, trim={0 0.5em 0 0},clip]{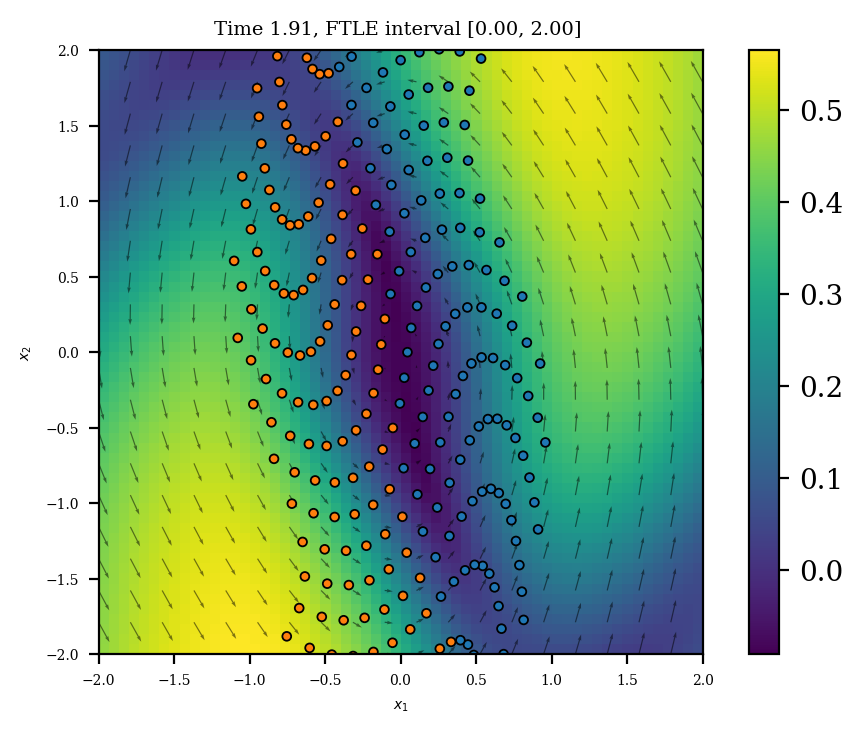}
	\end{subfigure}
		\begin{subfigure}{0.32\textwidth}
		\includegraphics[width = \textwidth, trim={0 0.5em 0 0},clip]{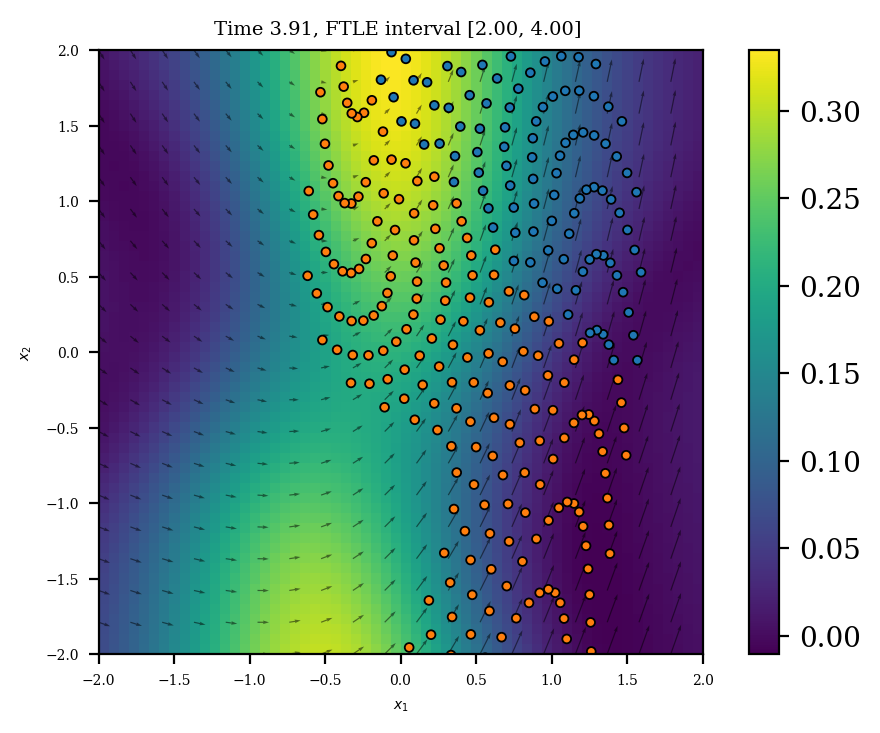}
	\end{subfigure}
		\\
		\begin{subfigure}{0.32\textwidth}
		\includegraphics[width = \textwidth, trim={0 0.5em 0 0},clip]{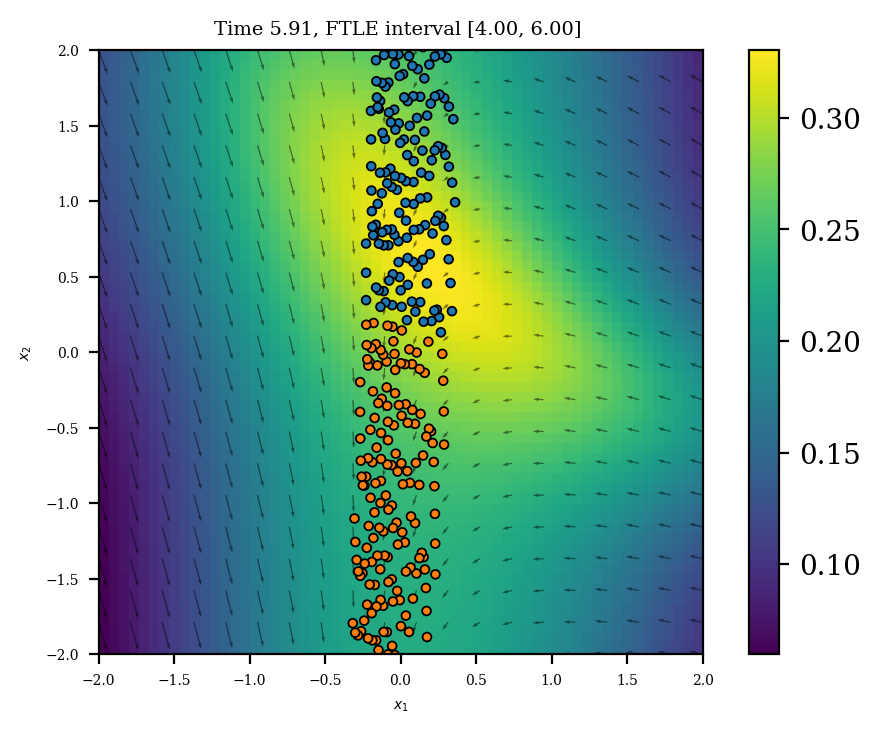}
	\end{subfigure}
		\begin{subfigure}{0.32\textwidth}
		\includegraphics[width = \textwidth, trim={0 0.5em 0 0},clip]{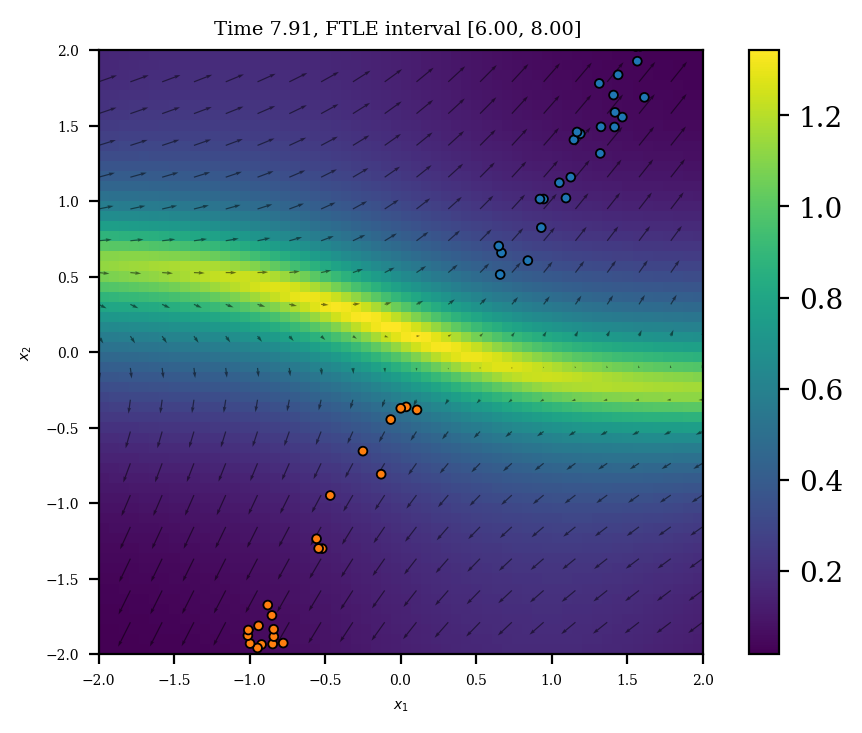}
	\end{subfigure}
	\begin{subfigure}{0.32\textwidth}
		\includegraphics[width = \textwidth, trim={0 0.5em 0 0},clip]{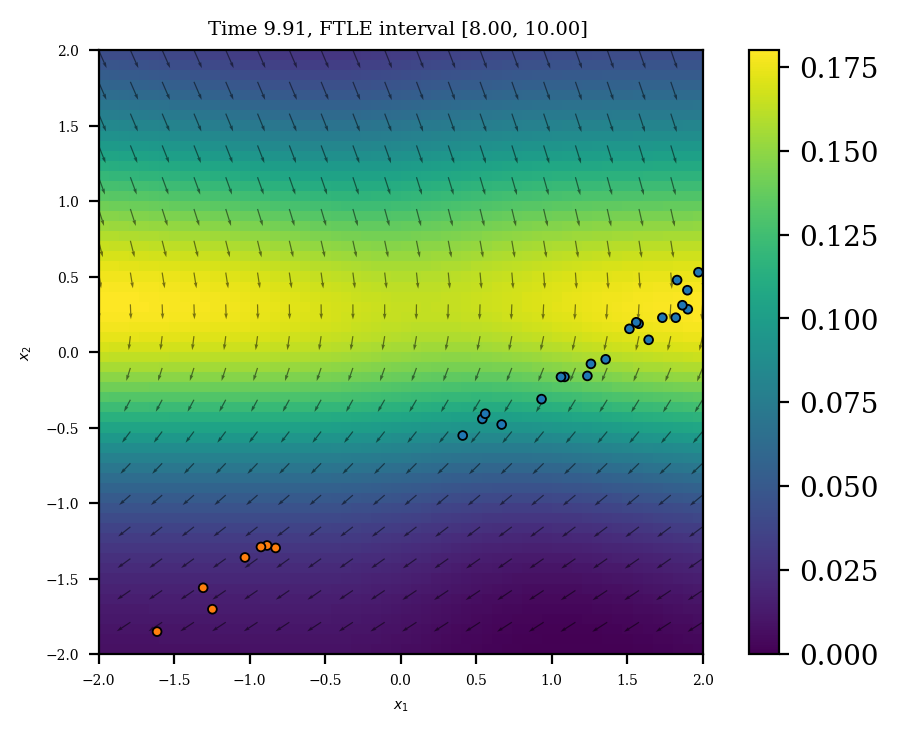}
	\end{subfigure}
	
	\caption{\small Plotted is the neural ODE evolution of Example~\ref{ex:5params} for a grid of initial points. The points are colored according to the models final prediction. The small arrows describe the piecewise-constant vector field that generate the dynamics. The background colors show the maximum FTLEs on the autonomous subinterval of the active vector field. \label{fig:subintervalFTLEsmoons}}
\end{figure}

\begin{figure}
	\centering 
	\begin{subfigure}{0.30\textwidth}
		\includegraphics[width=\textwidth, trim={0 0.5em 0 0},clip]{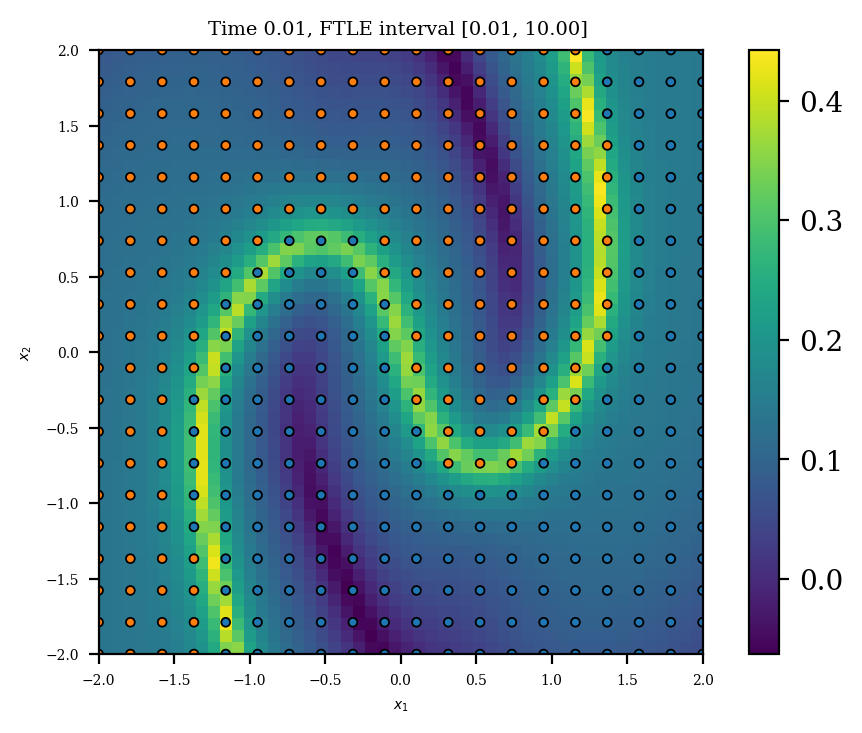}
	\end{subfigure}
	\begin{subfigure}{0.30\textwidth}
		\includegraphics[width = \textwidth, trim={0 0.5em 0 0},clip]{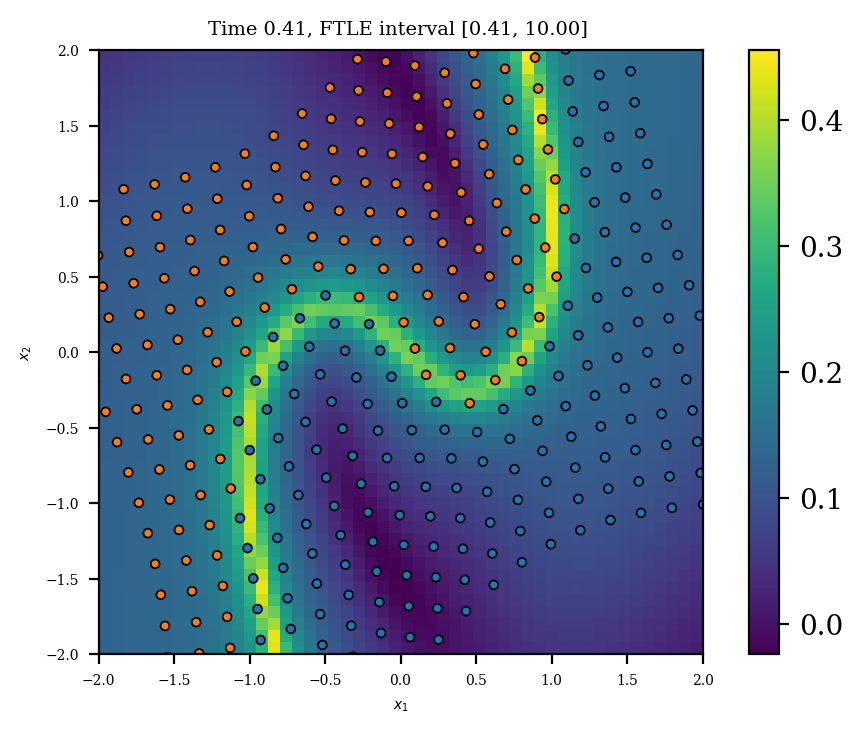}
	\end{subfigure}
	\begin{subfigure}{0.30\textwidth}
		\includegraphics[width = \textwidth, trim={0 0.5em 0 0},clip]{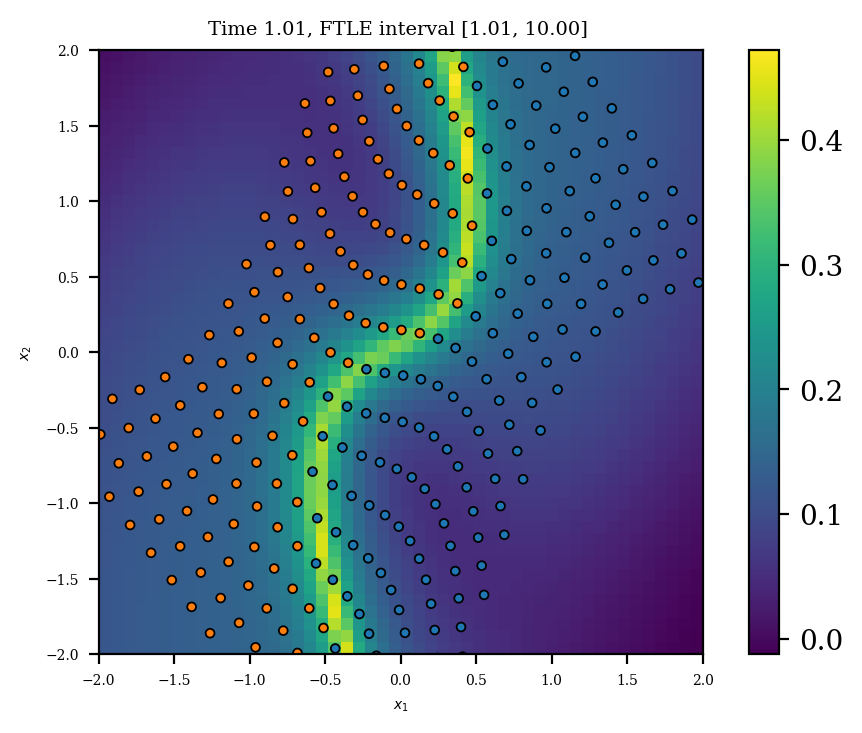}
	\end{subfigure}
	\\
	\begin{subfigure}{0.30\textwidth}
		\includegraphics[width = \textwidth, trim={0 0.5em 0 0},clip]{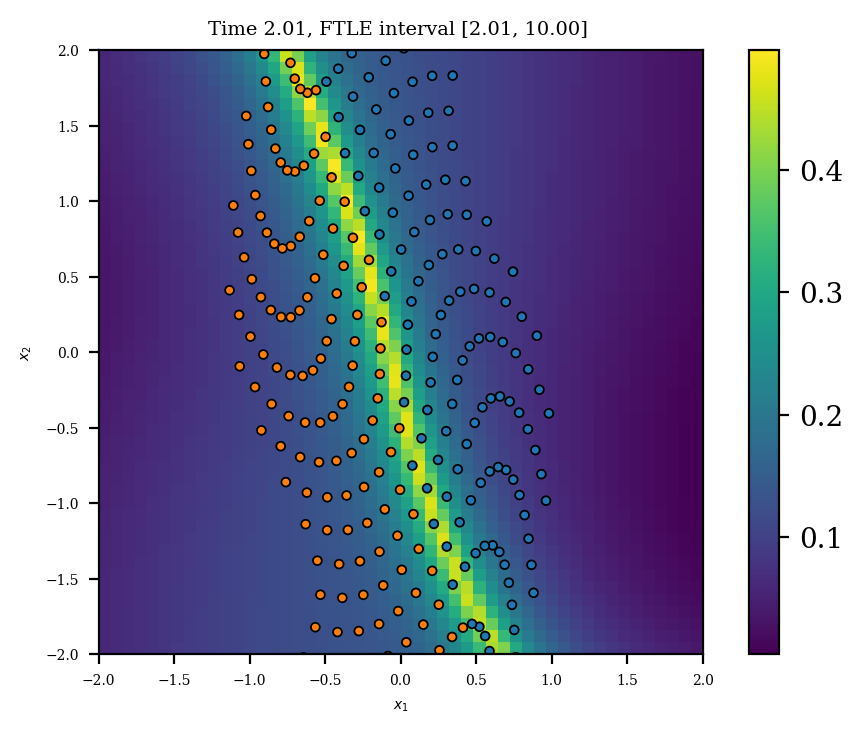}
	\end{subfigure}
	\begin{subfigure}{0.30\textwidth}
		\includegraphics[width = \textwidth, trim={0 0.5em 0 0},clip]{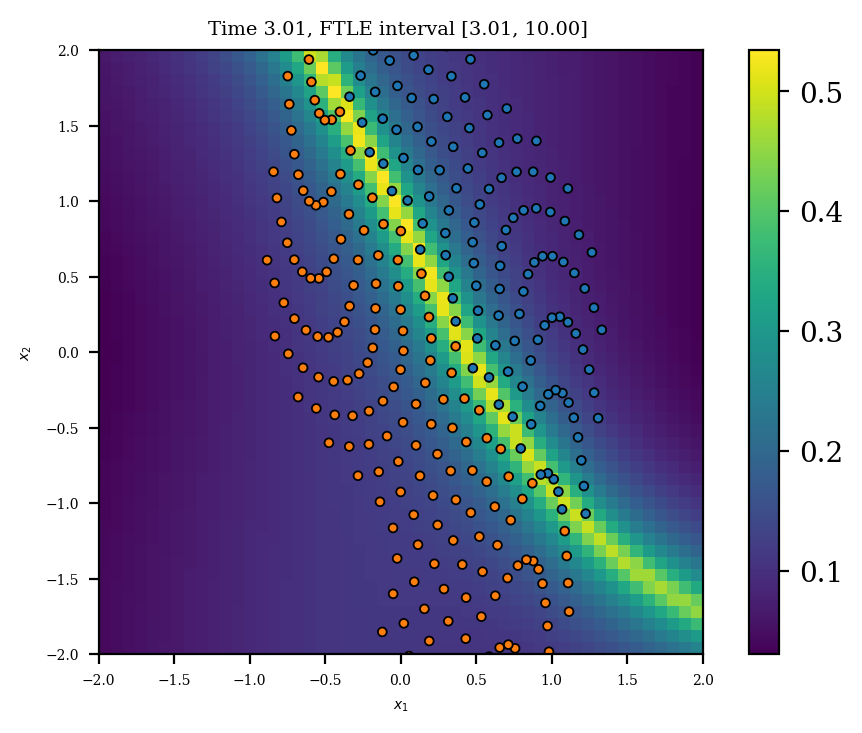}
	\end{subfigure}
	\begin{subfigure}{0.30\textwidth}
		\includegraphics[width = \textwidth, trim={0 0.5em 0 0},clip]{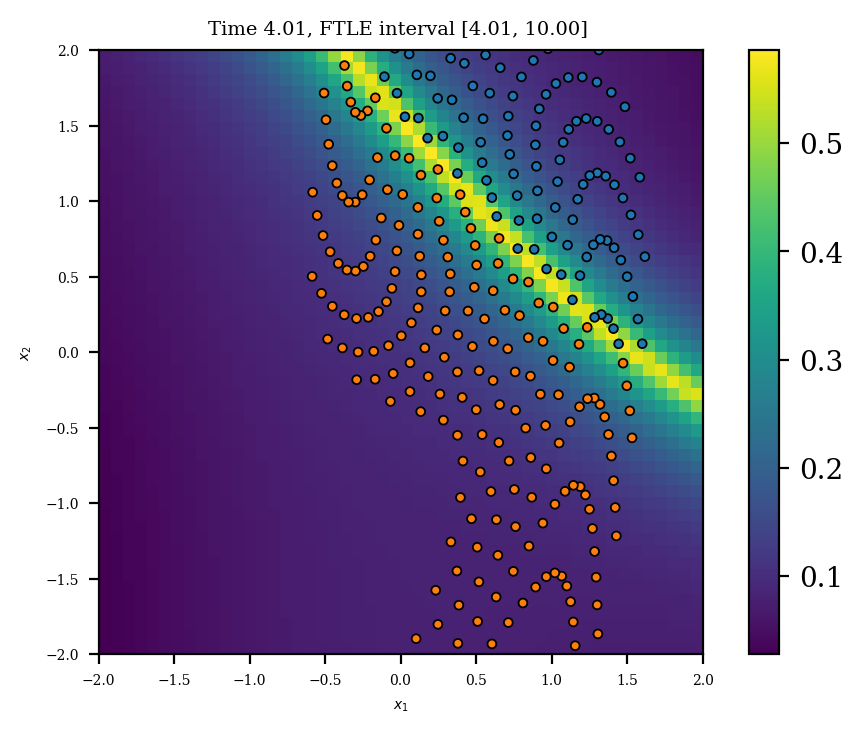}
	\end{subfigure}
		\\
	\begin{subfigure}{0.30\textwidth}
		\includegraphics[width = \textwidth, trim={0 0.5em 0 0},clip]{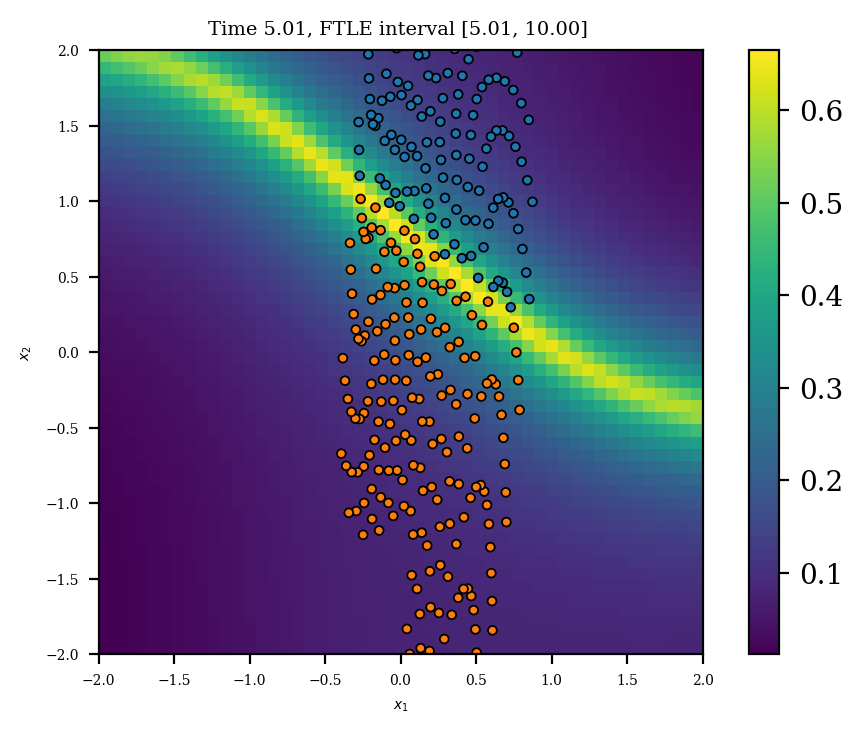}
	\end{subfigure}
	\begin{subfigure}{0.30\textwidth}
		\includegraphics[width = \textwidth, trim={0 0.5em 0 0},clip]{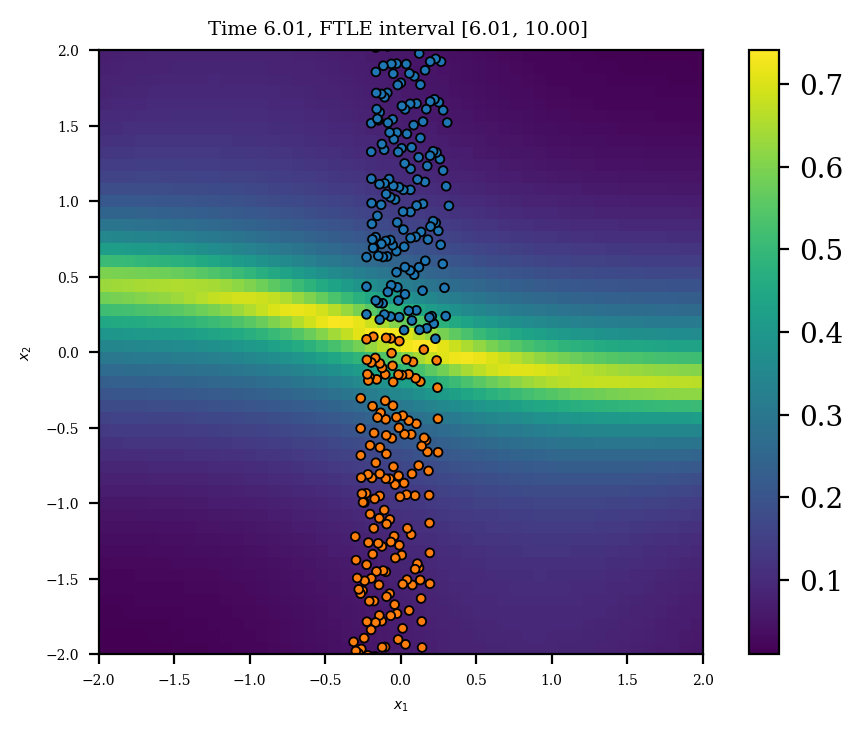}
	\end{subfigure}
	\begin{subfigure}{0.30\textwidth}
		\includegraphics[width = \textwidth, trim={0 0.5em 0 0},clip]{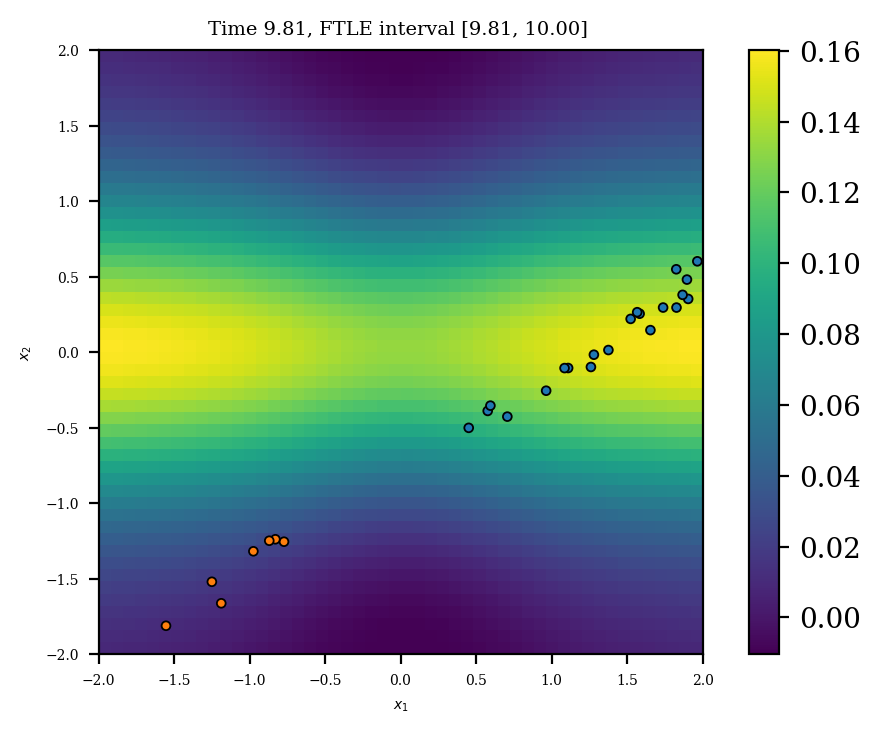}
	\end{subfigure}
	\caption{\small As in Figure~\ref{fig:subintervalFTLEsmoons}, plotted is the neural ODE evolution of Example~\ref{ex:5params}. However, here the background depicts the maximal FTLEs for the \emph{future time interval} which the points will traverse. \label{fig:shrinkingFTLE}}
\end{figure}

\setexample{2}
\begin{subexample}\label{ex:5params}
		We consider a non-autonomous neural ODE of type \eqref{eq:generalf} with five parameter layers, $K=5$, one vector field layer, $\ell = 1$, and all layer widths of constant size 2, $V_k = I$, $a_k = 0$ for $k = 1, \ldots 5$. The explicit vector field is given as
		\begin{equation*}
			f(t,x) = \sigma (W(t)x + b(t)), \quad t\in [0,10],
		\end{equation*}
		with 
		\begin{equation*}
			W(t) \equiv W_k \in \R^{2\times 2},\quad  b(t) \equiv b_k\in \R^2 \quad \text{for } t\in
            \textstyle\left[2(k-1), 2k\right), k=1,\ldots,5.
		\end{equation*}
Compared to Example~\ref{ex:1param}, at each time step this vector field has a simpler structure with fewer parameters active, operating as a sequence of five autonomous rather simple vector fields applied in series. However, its trajectories are less confined spatially and can revisit past regions in phase space due to the time-dependency. As in Example~\ref{ex:1param}, we discuss the dynamics of the trained nODE (with details on training described in Section~\ref{subsec:numerics}). Figure~\ref{subfig:moons:levelsets} shows that the trained model manages to achieve accurate classification on the moons dataset and manages to capture the desired generalization. Figure~\ref{subig:moons:ftles} again confirms that the decision boundary is captured by the FTLE ridge despite the structural differences. 
A notable difference appears when comparing the plots of the second Lyapunov exponents $\lambda_{\min}(x_0)$ for both examples. For Example~\ref{ex:1param} (Figure~\ref{subfig:moons:minftles1param}), the $\lambda_{\text{min}}(x_0)$-ridge corresponding to the decision boundary are local maxima of positive exponents, indicating divergence in all directions. Conversely, for Example~\ref{ex:5params} (Figure~\ref{subig:moons:minftles}), the $\lambda_{\min}(x_0)$-ridge are local minima of negative exponents. This suggests that trajectories near the decision boundary in the latter case behave similarly to trajectories near hyperbolic saddle points.

Our vector field is piecewise-autonomous, and Figure~\ref{fig:subintervalFTLEsmoons} provides the FTLEs of each autonomous subinterval that the inputs traverse. One observes that in the first subinterval, $[0,2]$, the most significant ``straightening'' of the nonlinear moon shapes is performed via spiraling dynamics. There, predominantly inputs away from the origin diverge, indicated by the brighter regions. The importance of the role of the later dynamics for the classification task is less evident. Overall, Figure~\ref{fig:subintervalFTLEsmoons} provides the nODE dynamics in detail, but it does not by itself yield a conclusive picture of the separation process. Based on these plots it is further not easy to observe how the dynamics of each subinterval combine to the FTLE ridge of the whole interval, plotted in Figure~\ref{subig:moons:ftles}.

Figure~\ref{fig:shrinkingFTLE} indicates the classification progress remarkably well and confirms the link of FTLE ridges to almost invariant sets described in Section~\ref{subsec:lcs}. Plotted are the maximum FTLEs for the \emph{future time interval} which the points will traverse, given as $\lambda_{\text{max}}((t,T],x_0)$ at current time $t\in (0,T)$. As time $t$ progresses, the ridge straightens from the initial nonlinear shape to an almost linear one. The plots show that the ridge associated to $\lambda_{\text{max}}((t,T],x_0)$ generates the boundaries of two almost invariant sets (as defined in \eqref{eq:tinvariant}) $S_{t,\text{orange}}$ and $S_{t,\text{blue}}$ corresponding to the model's classification. A similar observation can be made for Example~\ref{ex:1param}. The analogous plots are omitted as the ridge evolution is already contained in Figure~\ref{fig:growingFTLE1param} when read in reverse order. Indeed, then the nonlinear ridge again transforms into an almost linear one. 

Additionally, we again observe that the dynamics can be categorized roughly in two phases: An early phase, where the linearization of the nonlinear class boundary takes place, and a later phase, where the position of opposite class data points increases.\\
\end{subexample}

Let us summarize the observations of the above two examples:
\begin{itemize}
    \item The computed FTLE ridges coincide with the decision boundaries generated by the models. This is consistent for the two structurally different dynamics that realize the classification.
    \item The ridge associated to $\lambda_{\text{max}}((t,T],x_0)$ generates two almost invariant sets corresponding to the model's classification.
    \item The plot of the second FTLEs, $\lambda_{\text{min}}(x_0)$, also shows the FTLE ridge. However, for Example~\ref{ex:1param} both LE directions on the ridge are expanding, while Example~\ref{ex:5params} has one expanding and one contracting LE direction.
   \item For both the autonomous and the non-autonomous example, the dynamics can be split into the early dynamics which transforms the nonlinear decision boundary into a linear one, and the later dynamics which accentuate the spatial separation.
\end{itemize}

\section{Robust training via FTLE suppression}\label{sec:ftlereg}
We introduce adversarial robustness in the nODE setting and connect it to FTLE ridges. We then propose a modified training scheme that regularizes the loss function by suppressing large FTLEs, thereby robustifying the model against adversarial attacks. The neural ODE framework allows us to design this strategy from a dynamical systems perspective. Since a primary challenge lies in computational efficiency, we focus on the early dynamics, which we claim are both the dominant factor in adversarial vulnerability and reduce the computational load.

\subsection{Adversarial robustness}\label{subsec:adv}
Deep neural networks achieve high accuracy on test sets but often lack robustness on manipulated (or contaminated) inputs, called adversarial attacks. Such vulnerabilities severely restrict the real-world applications of DNNs. The problem is linked to the high dimensionality of applications and the constraint on computational resources in relation to the entire input space \cite{Good2015, madry2018}. Many strategies to mitigate the problem have been developed in recent years and can be grouped into two adversarial defense approaches. One approach involves adding more data to the training process, typically first-order estimates of maximal adversarial examples \cite{Good2015, madry2018}, while another adjusts the model structure, e.g., via regularization terms \cite{lyu15gradientreg, pt17ortho} or directly \cite{liu18randomlayer, schoenlieb25ortho}. We outline how adversarial vulnerability connects to (large) FTLEs and how to use that to improve the robustness of neural ODEs.

For the here considered classification tasks, adversarial robustness can be understood as follows. Given an input $x_0$, a model $\Psi$ aims to predict the class represented by a ground-truth label $y$, such that $\Psi(x_0) = y$.
An \emph{adversarial attack} with \emph{adversarial budget} $\eps>0$ takes the form $\tilde{x}_0= x_0 + \xi$ where $\|\xi\|\leq \eps$ in a predefined norm $\|\cdot\|$. The attack succeeds if $\Psi(\tilde{x}_0) \neq y$. Conversely, we call a model \emph{adversarially robust} at an input $x_0$ (with respect to $(\|\cdot\|,\eps)$-perturbations) if $\Psi(\tilde{x}_0) = y$ holds true for all attacks $\tilde{x}_0$ of budget $\eps$. In our framework, the model is a composition $\Psi(x_0) = L \Phi(T, x_0)$, where $\Phi(T,\cdot)$ denotes the time-$T$ flow defined in \eqref{eq:flow} and $L$ is an affine layer defined in \eqref{eq:loss}. It follows that adversarial vulnerability is tied to high local input sensitivity of the input output map $\Psi$. The linearized sensitivity of an input $x_0$ in direction $\xi\in \R^d$ is given as
\begin{equation}
\lim_{\eps \to 0} \frac{\Psi(t, x_0 + \eps \xi)  - \Psi(t, x_0)}{\eps} = D_xL \,D_{x_0}\Phi(t,x_0) \xi.
\end{equation}
Dynamically, we can neglect the role of the final layer $L$ here and focus on the sensitivity analysis of the nODE flow.
Choosing the direction of the maximum LE, $\xi = \xi_{\max}$, we obtain (cf.\ \eqref{eq:deflmax} and \eqref{eq:cauchygreen})
\begin{equation}
    \lambda_{\max}(t,x_0) = \frac{1}{t}\lim_{\eps \to 0} \ln\left(  \frac{\|\Phi(t,x_0 + \eps \xi_{\max}) - \Phi(t, x_0)\|_2}{\eps} \right).
\end{equation}
Consequentially, the maximum LE for final time $T$ provides a direct measure of the (first order maximum) input-output sensitivity of the neural ODE and as such a local indicator of adversarial vulnerability in Euclidean norm.

\begin{figure}
\includegraphics[scale = 0.7]{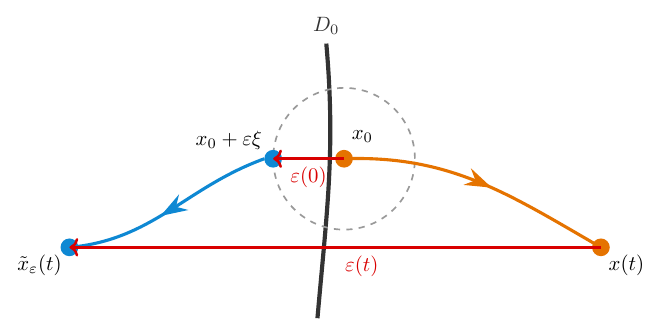}
    \centering
    \caption{\small Given is an input $x_0$ close to the decision boundary $D_0$ (black curve). The input perturbation $x_0 + \eps \xi$ represents a successful adversarial attack which leads to an input with opposing class prediction. Dynamically, the perturbation distance grows with time, $\eps(0) \ll \eps(t)$, resulting in a large positive Lyapunov exponent at the input $x_0$.}
    \label{fig:decisionboundary}
\end{figure}

High input sensitivity is typically observed near the model's decision boundaries, $D_0$ defined in \eqref{eq:Deps}. At these boundaries large maximum FTLEs are expected, as the classification is realized by pushing inputs of opposing class predictions into opposing directions, see Figure~\ref{fig:decisionboundary}. The consequence is that the ``sharper'' the decision boundaries, the larger the maximum FTLEs nearby. 

Among other possible approaches, our strategy is to robustify models by ``blurring'' or widening the decision boundary and the associated FTLE ridges. This reduces local input sensitivity and creates a buffer of low prediction confidence between the classes, where adversarial attacks may be more readily detected.

\subsection{Modified loss function}\label{subsec:ftlereg}
As discussed in Section~\ref{subsec:adv}, large FTLEs introduce input perturbation sensitivity into our dynamics which can have the unwanted effect of adversarial vulnerability. The goal here is to improve adversarial robustness for neural ODEs of type \eqref{eq:node}--\eqref{eq:generalf} by suppressing large FTLEs via a modified loss function. More specifically we aim to reduce high-confidence mistakes realized via a ``blurring'' or broadening of the decision boundary, see Figure~\ref{fig:ftlereg1param}.

We introduce the modified loss function with added \emph{regularization term} of strength $\gamma \geq 0$ as
		\begin{equation}\label{eq:totalloss}
		\mathcal{L}_\gamma({\mathbf {x_0,y}}) := \mathcal{L}({\mathbf{ x_0,y}}) + \gamma \mathcal{R}(\theta,{\mathbf{x_0}} ),
		\end{equation}
with $\mathcal{L}$ from \eqref{eq:loss} and the \emph{FTLE regularization term}
\begin{equation}\label{eq:ftlereg}
	\mathcal{R}(\theta,\mathbf{x_0}): = \frac{1}{N_{\text{data}}}\sum_{i = 1}^{N_\text{data}} \max \{\lambda_{\max}([0,T_1],x_{0,i}), \delta\}
\end{equation}
        with the tuneable hyperparameters of the threshold $\delta > 0$ and the time horizon $T_1 \in (0,T]$. The threshold parameter $\delta$ sets the cutoff for suppressing positive Lyapunov exponents. The time $T_1$ defines the duration of the neural ODE integration considered for the FTLE regularization. We refer to $\gamma = 0$ as the standard training/loss function case.

        As mentioned in Section~\ref{subsec:adv}, a large variety of regularization strategies exist. Here, we focus solely on FTLE suppression to isolate and analyze its effect on the dynamics. Closely related regularization strategies based on Lyapunov exponents have been suggested in \cite{engelken23, vockmeisel25}. Identity \eqref{eq:deflmax} further reveals a close relation of FTLE regularization to established gradient regularization strategies \cite{doubletrillos, lecun92, lyu15gradientreg}.

        For the two-dimensional setting analyzed here, the strategy was sufficient to \emph{only} regularize the first Lyapunov exponent $\lambda_{\max}$, rather than including also the second $\lambda_{\min}$. For higher dimensional problems, we expect the inclusion of a small finite number (compared to the input dimension) of leading LEs to be a reasonable choice. Other FTLE based regulariziation terms, such as $\lambda_{\max}([0,T_1],x_0)^2$ \cite{engelken23}, are worth considering but here we aimed to show the isolated effect of reducing positive FTLEs.
 
\subsection{FTLE suppression in two-dimensions}\label{subsec:2dftlereg}

We continue with Example~\ref{ex:1param} and Example~\ref{ex:5params} but here perform a modified training with suppression of the Lyapunov exponents as described in \eqref{eq:totalloss}. To clearly observe the effect of the modified training, we compare two models per example with identical parameter initialization. One trained with standard loss, $\gamma = 0$ in \eqref{eq:ftlereg}, and one with an activated FTLE suppression term \eqref{eq:ftlereg}.
\\

\begin{figure}
	\centering 
	\begin{subfigure}{\textwidth}	\centering 
		\includegraphics[width = 0.32\textwidth, trim={0 2.3em 0 0},clip]{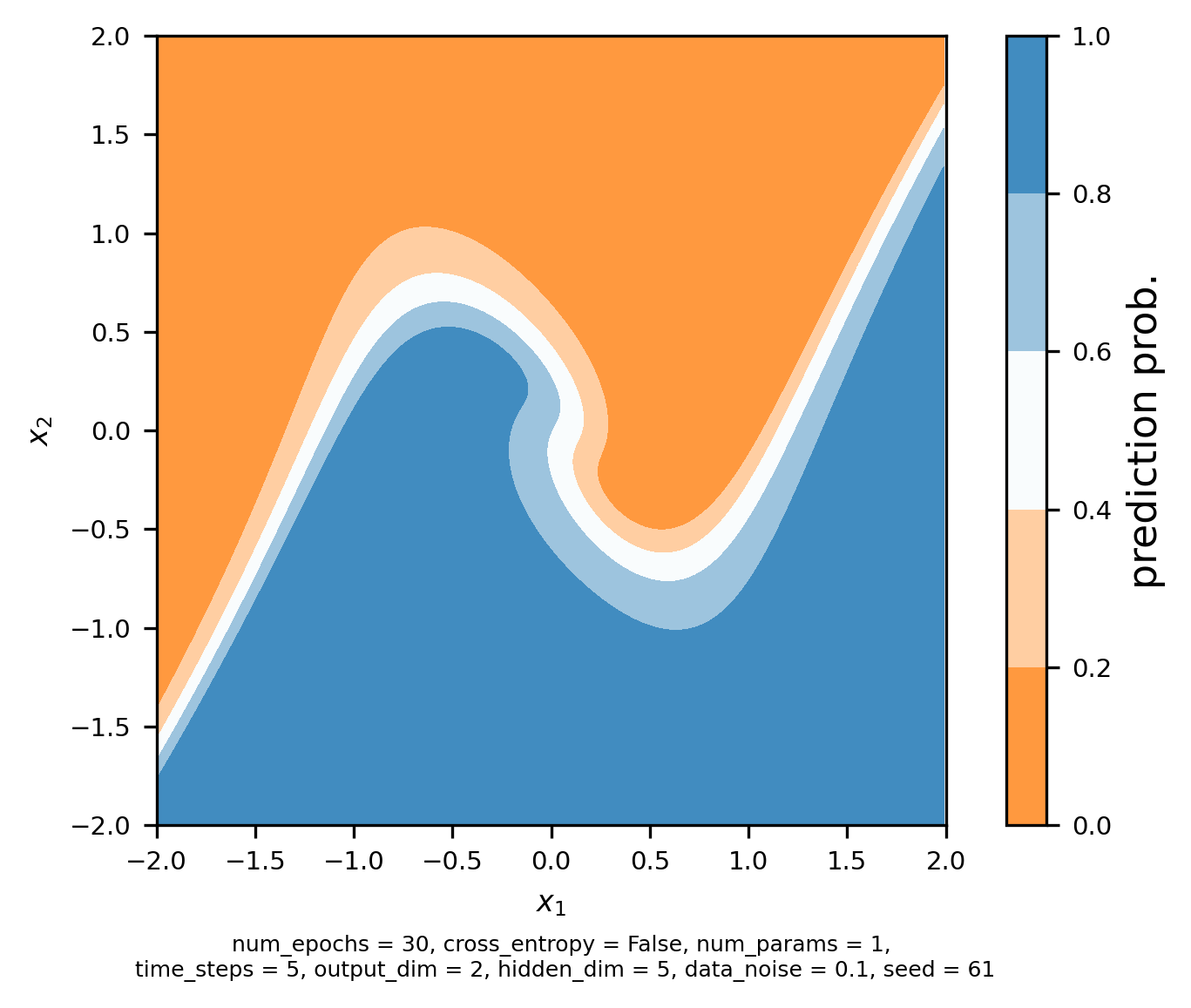}
		\includegraphics[width = 0.32\textwidth, trim={0 0.em 0 0},clip]{./plots/regularization/1param/MLE_max_control}	\includegraphics[width = 0.26\textwidth, trim={0 0.em 0 0},clip]{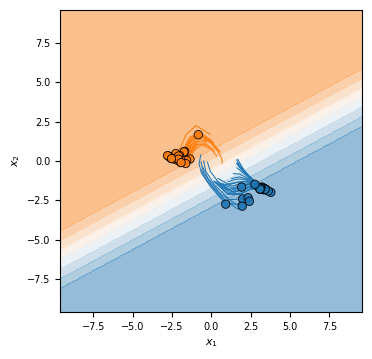}
		\caption{\small Model without FTLE suppression \label{subig:mftlenoreg}}
	\end{subfigure}\\
	\begin{subfigure}{\textwidth}	\centering 
		\includegraphics[width = 0.32\textwidth, trim={0 2.3em 0 0},clip]{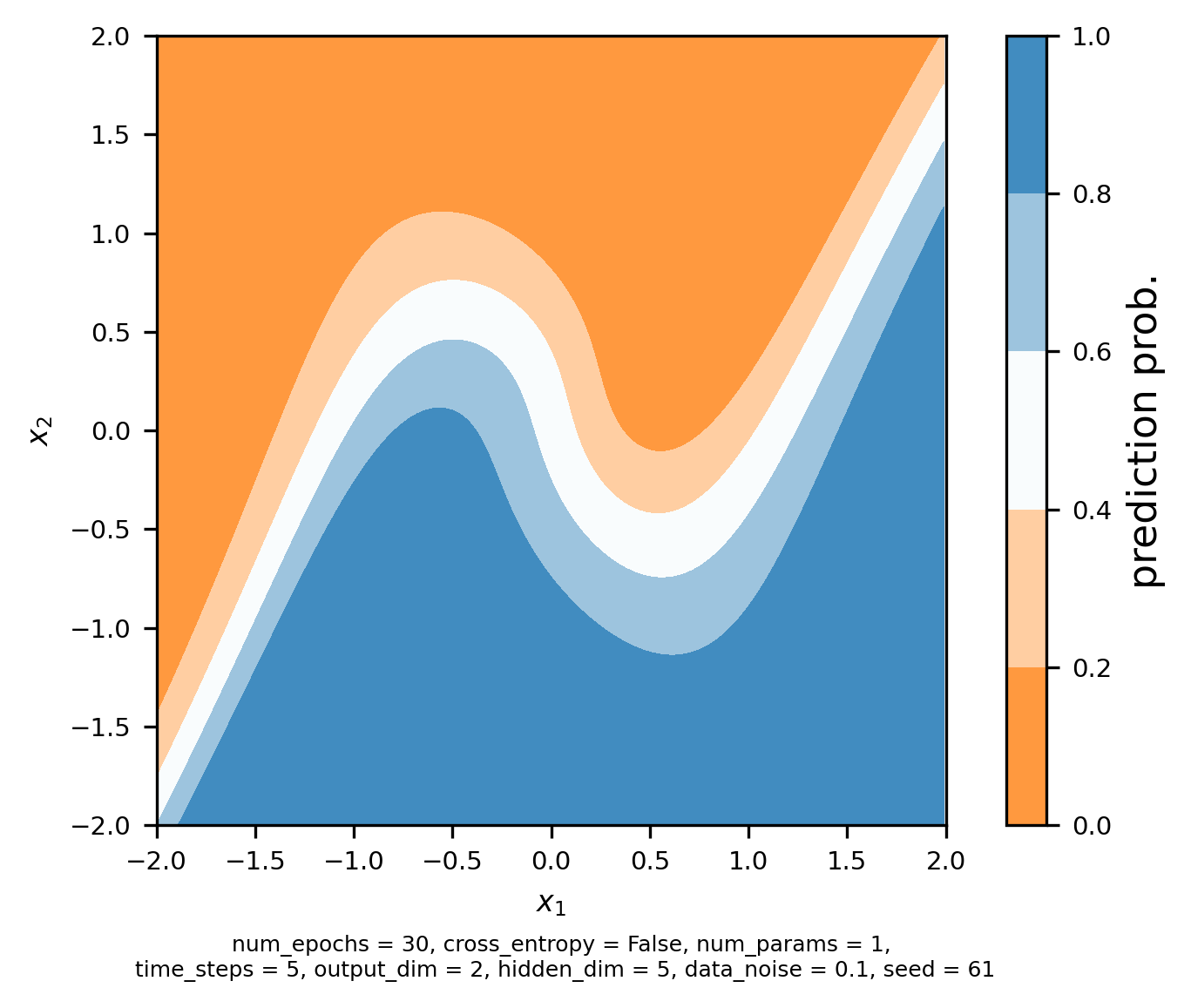}
		\includegraphics[width = 0.32\textwidth, trim={0 0.em 0 0},clip]{./plots/regularization/1param/MLE_max}
		\includegraphics[width = 0.26\textwidth, trim={0 0.em 0 0},clip]{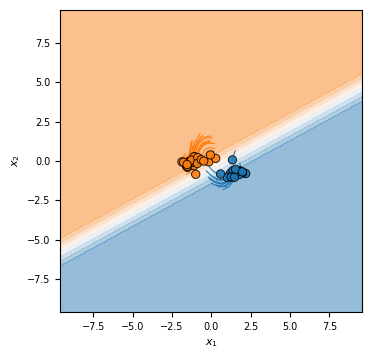}
		\caption{\small Model with FTLE suppression on $[0,2]$ \label{subig:trajreg}}
	\end{subfigure}
	\caption{\small Comparison of Example~\ref{ex:1paramcont} trained without and with active FTLE suppression\label{fig:ftlereg1param}}
\end{figure}	

\setexample{1}
\setcounter{subexample}{1} 
\begin{subexample}[Modification of Example~\ref{ex:1param}]\label{ex:1paramcont}
We set $\delta = 0.05$, regularization strength $\gamma = 2$ and $T_1>0$ is chosen according to the subcases specified below.
\begin{itemize}
\item \emph{FTLE regularization on the whole interval, $T_1 = T$:} Compared to standard training this decreases the maximal Lyapunov exponents $\lambda_{\max}([0,10],x_0)$ almost everywhere equally. The dynamics change, however not its general structure, resulting simply in shorter trajectories. The classification level sets are almost identical to the ones trained without the regularization term, as $L$ is able to compensate the changed dynamics. This confirms that we can have very similar classification results for rather different FTLE plots. The FTLE plot gives us an extra layer of information on the dynamics on top of the classification level sets and can serve as an additional criterion when selecting models.

\item \emph{FTLE regularization on the early dynamics, $T_1 = 2$, see Figure~\ref{fig:ftlereg1param}:} The time horizon choice $T_1 = 2$ is informed by the insight of Experiment~\ref{ex:1param}, specifically the discussions on Figure~\ref{fig:growingFTLE1param}. It reduces large input divergence of the fast dynamics, which has the most impact on separating the classes. But it still allows enough expressivity of the dynamics to carry out the classification process. The effect is that the decision boundary $D_\eps$ with $\eps = 0.1$ (defined in \eqref{eq:Deps}) widens in the regularized case, as visible when comparing the white strips in the first column. The second column shows that the $\lambda_{\max}([0,10],x_0)$-ridge around the origin disappeared. This reduces local input sensitivity while the model still separates dataset inputs successfully. It maintains class prediction accuracy and improves adversarial robustness.
\end{itemize}

 We conclude that for our example with a vector field of form \eqref{eq:vectorfieldex1}, an FTLE suppression targeting only the input's early dynamics is effective in avoiding high-confidence mistakes around the decision boundary. Compared to active suppression for all times, a suppression restricted to the early dynamics also has the advantage of being computationally less demanding, discussed in more detail in Section~\ref{subsec:ftleregeff} below.
\end{subexample}

\begin{figure}
	\centering 
	\begin{subfigure}{\textwidth}	\centering 
		\includegraphics[width = 0.325\textwidth, trim={0 2.3em 0 0},clip]{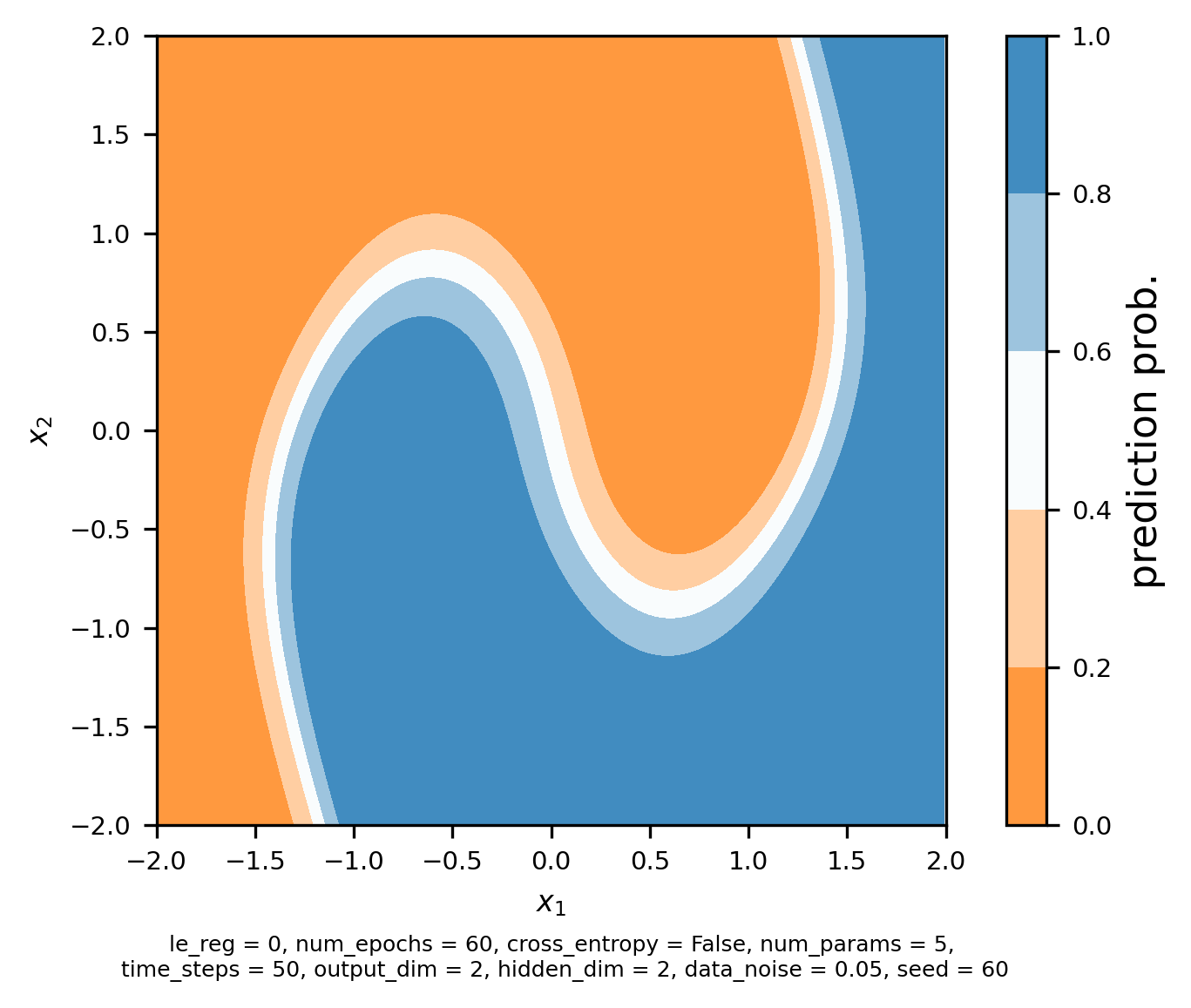}
		\includegraphics[width = 0.32\textwidth, trim={0 0.em 0 1.8em},clip]{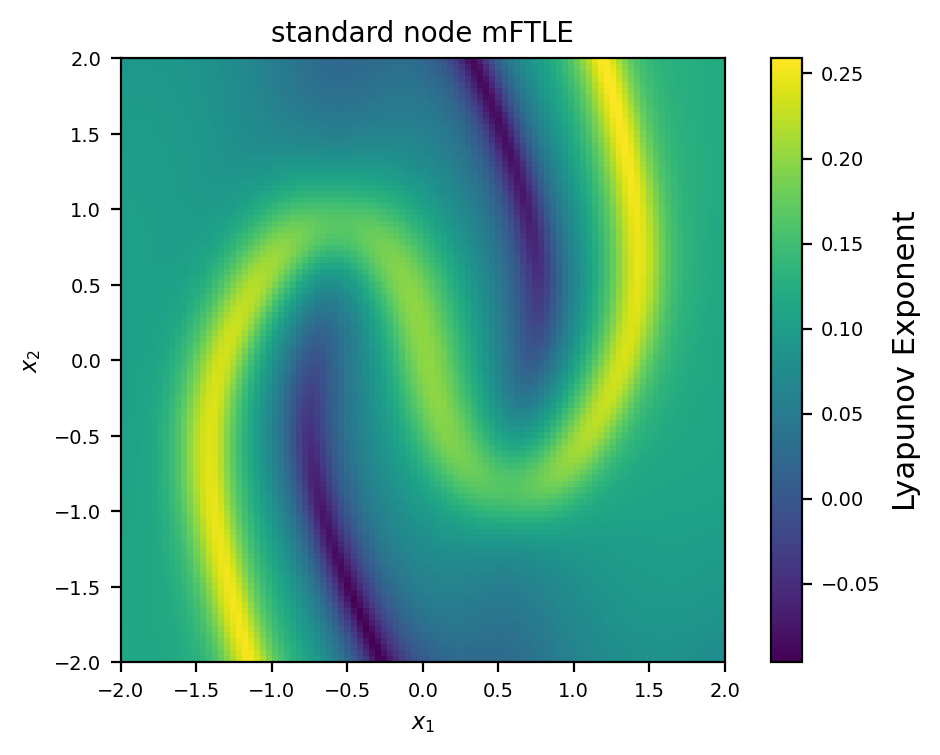}	\includegraphics[width = 0.25\textwidth, trim={0 0.em 0 0},clip]{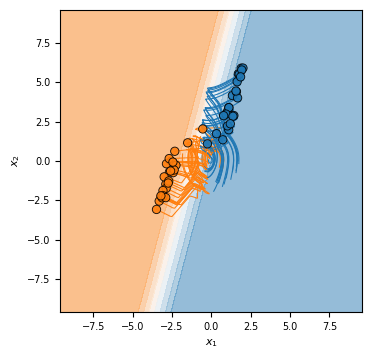}
		\caption{\small Model without FTLE suppression \label{subig:5params:mftlenoreg}}
	\end{subfigure}\\
	\begin{subfigure}{\textwidth}	\centering 
		\includegraphics[width = 0.325\textwidth, trim={0 2.3em 0 0},clip]{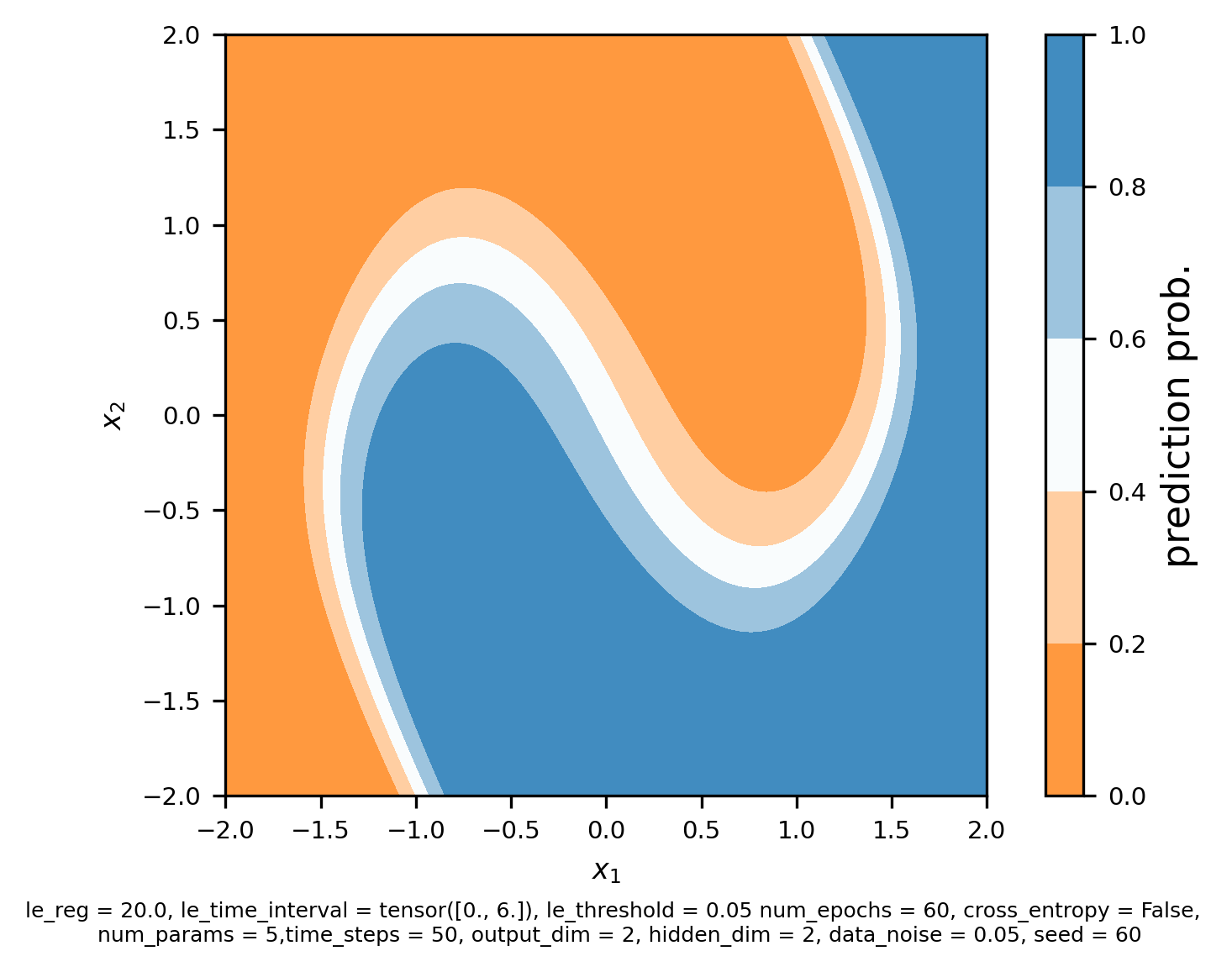}
		\includegraphics[width = 0.31\textwidth, trim={0 0.em 0 1.8em},clip]{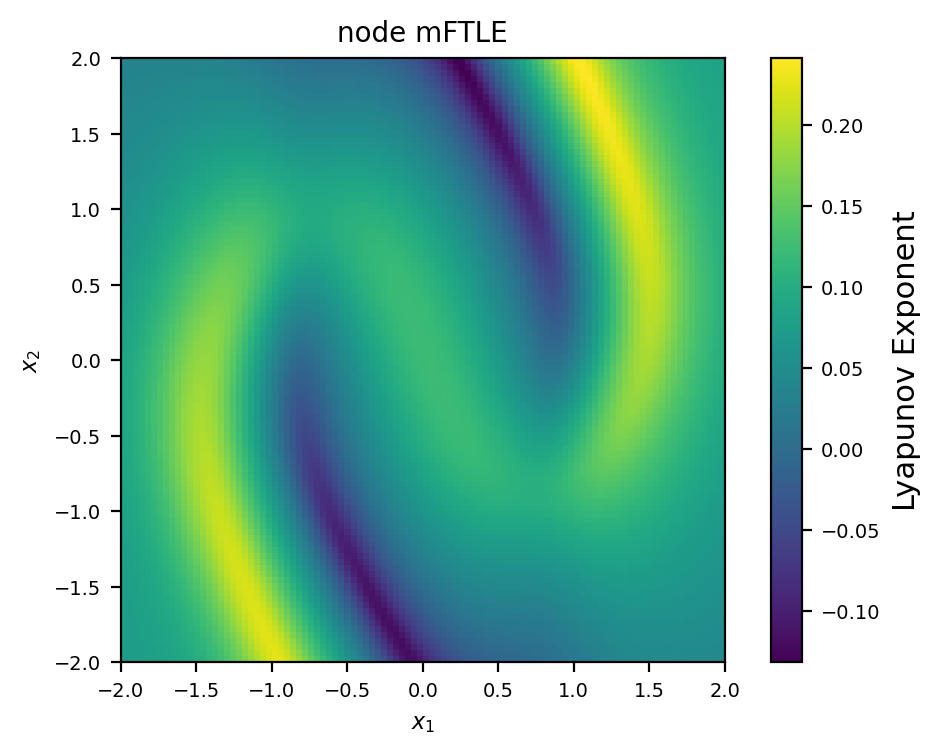}
		\includegraphics[width = 0.25\textwidth, trim={0 0.em 0 0},clip]{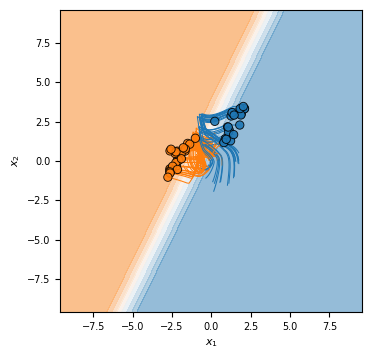}
		\caption{\small Model with FTLE suppression on $[0,6]$. \label{subig:5params:trajreg}}
	\end{subfigure}
	\caption{\small Comparison of Example~\ref{ex:5paramscont} without and with active FTLE suppression. \label{fig:ftlereg5params}}
\end{figure}	

\setexample{2}
\setcounter{subexample}{1} 

~
\begin{subexample}[Modification of Example~\ref{ex:5params}]\label{ex:5paramscont}
We set $\delta = 0.05$, $\gamma = 20$ and $T_1>0$ is chosen as specified below.
\begin{itemize}
    \item \emph{FTLE regularization on the whole interval, $T_1 = T = 10$}: As above in Example~\ref{ex:1paramcont}, the case $T_1 = 10$ does not result in significant changes of the decision boundary or structure of the dynamics and merely slows down the training process due to the high computational load and suppressed dynamical expressivity. 
\item \emph{FTLE regularization on the early dynamics, $T_1 = 6$}:
The choice $T_1 = 6$ regulates the FTLEs corresponding to the first two parameter layers. The modified training is shown in Figure~\ref{fig:ftlereg5params}. Even though the dynamics are structured fundamentally differently to Example~\ref{ex:1paramcont}, the achieved effect of FTLE suppression is very similar: Compared to the standard training, the decision boundary is wider which helps to mitigate high-confidence misclassifications. Plotting $\lambda_{\max}([0,10],x_0)$ shows that the ridge in the center of the image is blurred, as seen in the middle column of Figure~\ref{fig:ftlereg5params}. The right column of Figure~\ref{fig:ftlereg5params} shows that the suppressed dynamics move the final trajectory positions $\Phi(T,x_0)$ closer together. This is happening predominantly in the direction perpendicular to the level set boundaries generated by the final affine-linear layer $L$. Figure~\ref{fig:subFTLEsupp} plots the input-output dynamics for the first autonomous subintervals. It shows that despite identical initialization, the standard and modified training can lead to rather different dynamics that realize the classification. This is particularly evident when comparing the plots of the second column of Figure~\ref{fig:subFTLEsupp} where the general flow direction and local FTLE ridge differ significantly.
\end{itemize}

 \begin{figure}
 	\centering 
 	\begin{subfigure}{\textwidth}
 		\includegraphics[width = 0.3\textwidth, trim={0 0.5em 0 0},clip]{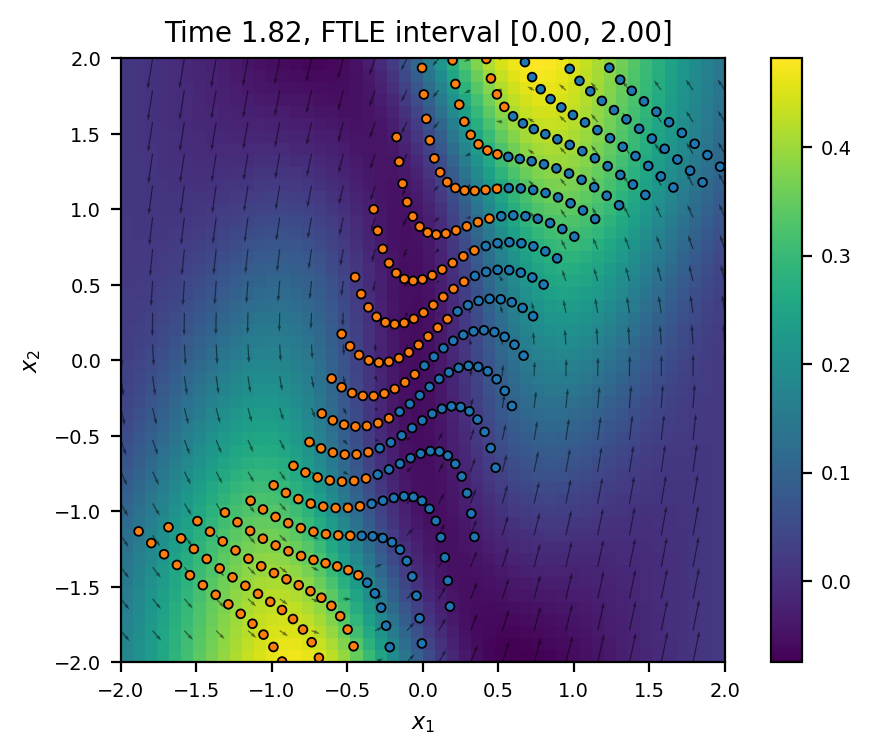}
 		\includegraphics[width = 0.3\textwidth, trim={0 0.5em 0 0},clip]{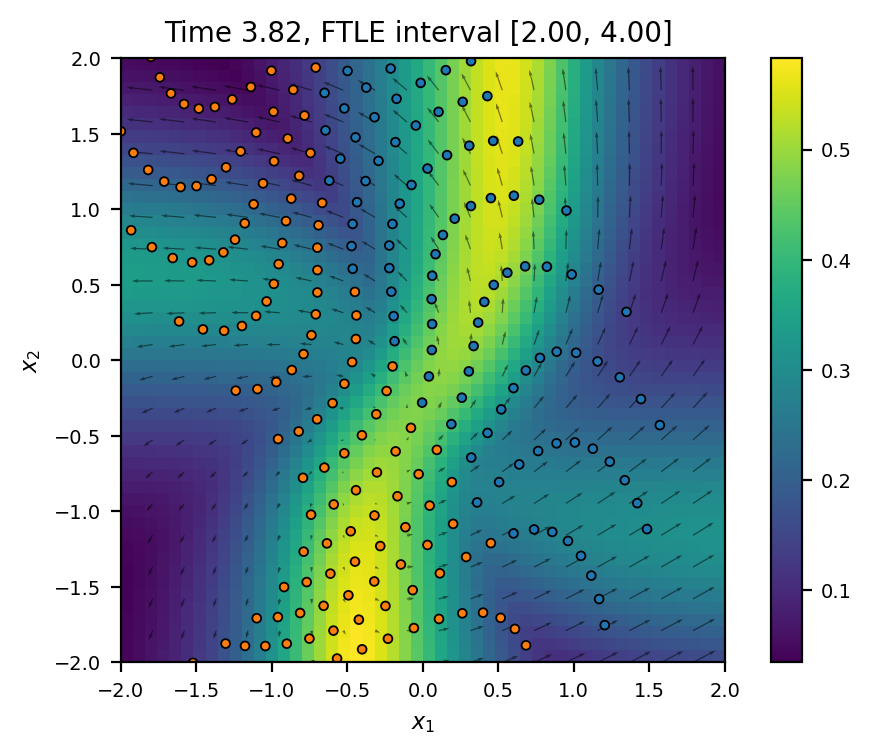}
 		\includegraphics[width = 0.3\textwidth, trim={0 0.5em 0 0},clip]{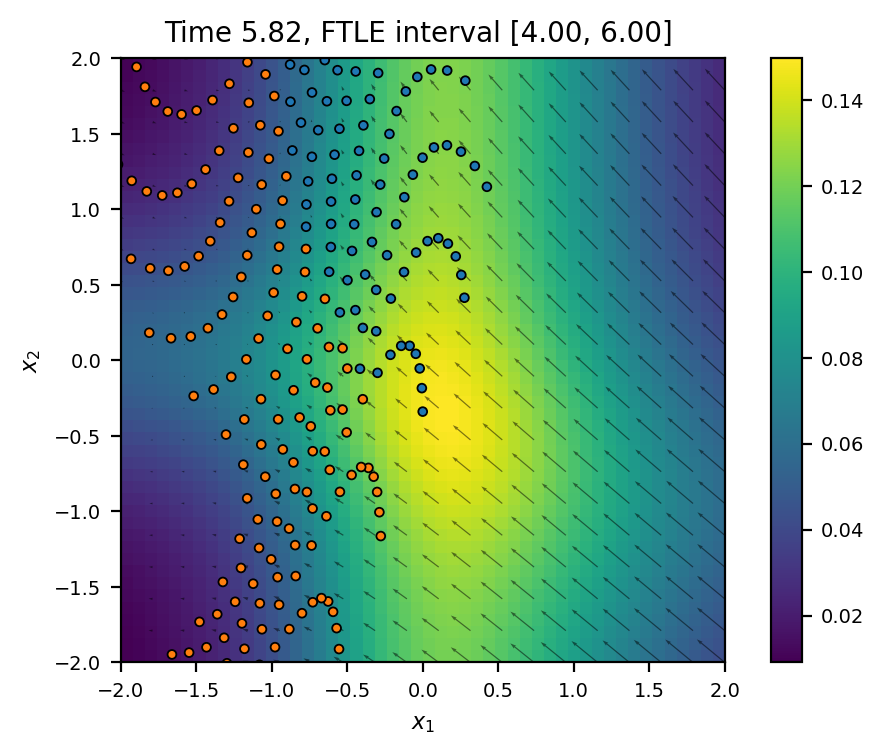}
 		\caption{\small  Model without FTLE suppression.}
 	\end{subfigure}
 	\begin{subfigure}{\textwidth}
 		\includegraphics[width = 0.3\textwidth, trim={0 0.5em 0 0},clip]{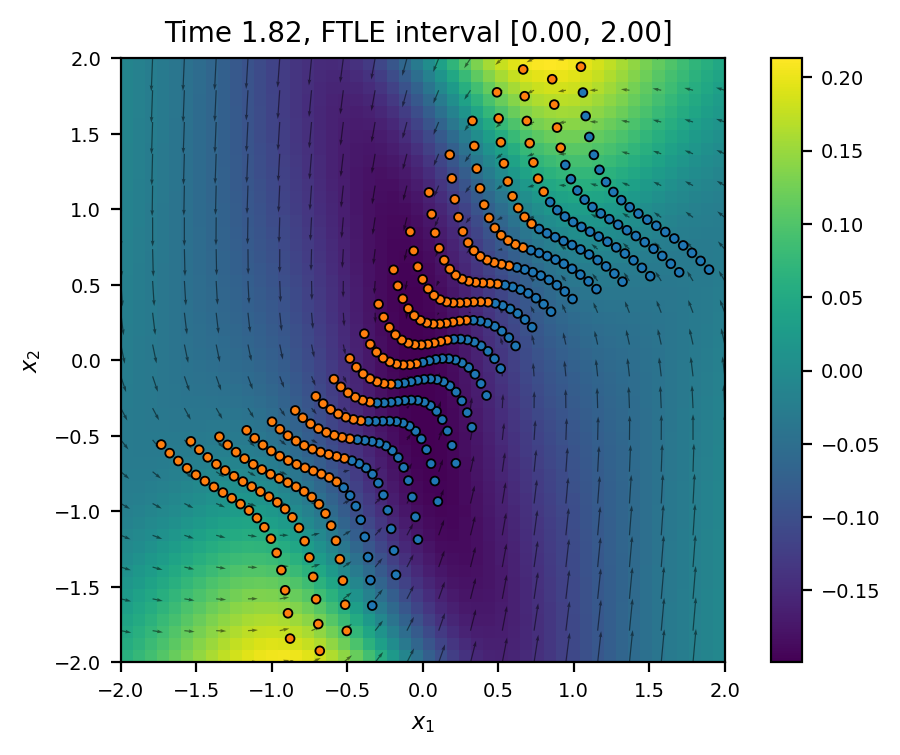}
 		\includegraphics[width = 0.3\textwidth, trim={0 0.5em 0 0},clip]{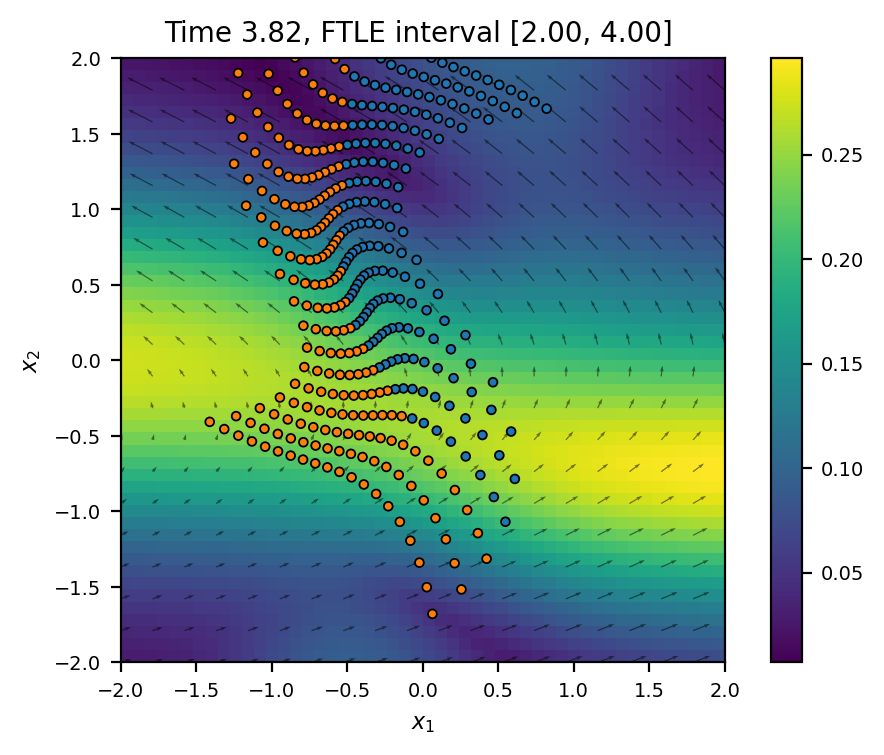}
 		\includegraphics[width = 0.3\textwidth, trim={0 0.5em 0 0},clip]{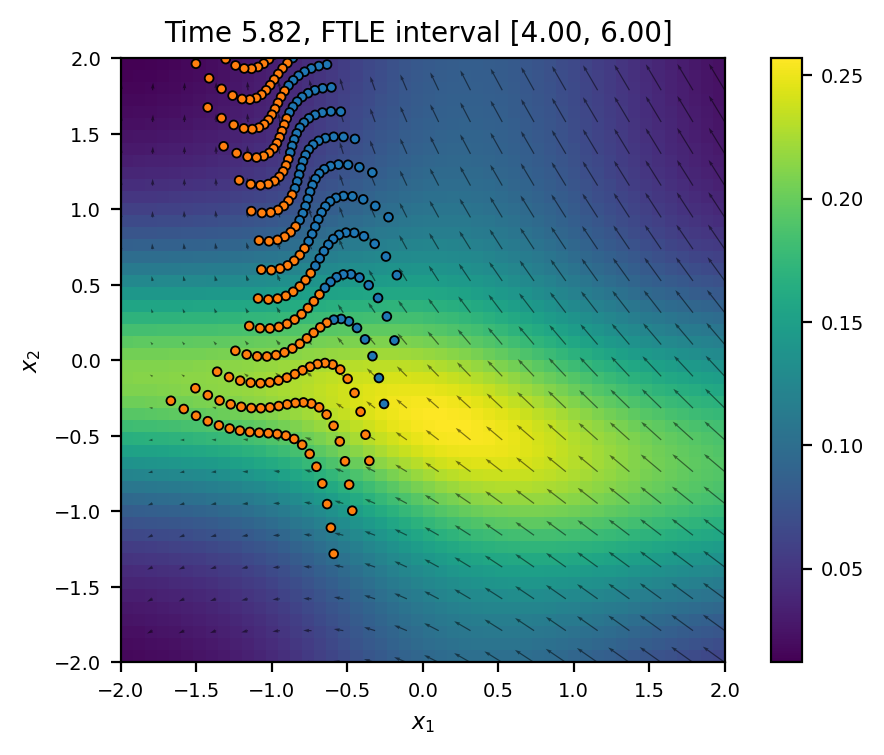}
 			\caption{\small  Model with FTLE suppression.}
 	\end{subfigure}
 	\caption{\small As in Figure~\ref{fig:subintervalFTLEsmoons}, we plot the neural ODE evolution of Example~\ref{ex:5paramscont}. \label{fig:subFTLEsupp}}
 	\end{figure}

We conclude Example~\ref{ex:5paramscont} by observing that a similar strategy as for Example~\ref{ex:1paramcont} to reduce high FTLEs on the early dynamics also is effective in this setting.
\end{subexample}

\begin{remark}
	We point out again that the final affine linear layer $L$ of \eqref{eq:loss} is not included in the regularization procedure. However, as $L$ has learnable parameters it may adjust to the suppressed Lyapunov exponents. We observe this on the right column of Figure~\ref{fig:ftlereg1param}, which shows that the $L$-generated level sets shrink closer together as the dynamics are less expansive. Therefore, it is necessary to view reduced FTLEs always in relation to the final layer. This also contributes to the explanation of why a regularization on the whole time interval typically has little effect on the level sets. It just slows down all dynamics while $L$ adjusts accordingly without changing the structure of the dynamics.
\end{remark}

\subsection{Efficient implementation of FTLE suppression}\label{subsec:ftleregeff}
Computation of Lyapunov exponents has been optimized for many years \cite{DRV97} and shall not be the focus here. One can readily apply available algorithms to our setting, such as the $QR$ decomposition based approach \cite{engelken2020lyapunov, prabhu2018}. Our parameter training is based on (stochastic) gradient descent for the loss \eqref{eq:totalloss} with regularizer $\mathcal{R}$ defined in \eqref{eq:ftlereg}. 
The update rule for a parameter $\theta_n$ in the parameter space based on a single input $x_0$ is computed as
\begin{align*}
\theta_{n+1}
&= \theta_n + s_n D_\theta \mathcal{L}_\gamma(\theta_n, x_0)\\
&= \theta_n + s_n D_\theta \mathcal{L}(\theta_n, x_0)
+ s_n \gamma \D_\theta \max \{\lambda_{\max}([0,T_1];x_0), \delta\},
& 
\end{align*}
where $s_n\!\ge 0$ is the learning rate.
This means each step requires the computation of two gradients: that of the standard square-loss gradient and that of the FTLE regularization term. The square-loss gradient can be realized in the ``naive'' way of backpropagating (i.e., applying the chain rule of differentiation) through the operations of the nODE integrator or via adjoint methods for which we refer to \cite{CRBD18, gholami2019}. The second term based on the regularization deserves a thorough discussion. To compute $D_\theta \mathcal{R}(\theta, x_0)$ we require:
\begin{itemize}
	\item An evaluation of $\mathcal{R}(\theta,x_0)$: Per definition of $\lambda_{\max}([0,T_1];x_0)$, this requires  the norm of the Jacobian, $\|\D_{x_0}\Phi(T_1,x_0)\|_2$, see Section~\ref{sec:ftledef}. 
	\item The derivative $D_\theta D_{x_0}\Phi(T_1,x_0)$: To obtain the parameter gradient of the regularizer another differentiation is required, now with respect to $\theta$.
\end{itemize}
Practically, the above steps correspond to a ``\emph{double backpropagation}'' (cf.\ \cite{lecun92,doubletrillos}), but only of the reduced flow $\Phi(T_1,\cdot)$: One pass from $T_1$ to $0$ to get the input sensitivity, and a second pass through the input sensitivity's dependence on parameter changes.  Also note that, in the case of a layered vector field $\ell > 0$ (such as Example~\ref{ex:1paramcont}), the backpropagation passes through the vector field layers in each discretized time step.
 It follows that the FTLE regularization term gradients \emph{require a computational load that is quadratically dependent on} $T_1$. It is therefore at the same time advantageous to select a small $T_1$ to improve adversarial robustness but also to significantly reduce the computational load of the regularization (compared to the regularization of the full dynamics).
 
The dynamical viewpoint permits another flexibility:
It is possible to use different time discretization schemes to optimize the square loss term and the regularization term. For example, a coarser discretization of the flow might be sufficient to compute the Lyapunov exponent regularization term. At the same time finer discretization of the square loss may improve the classification accuracy.

The here described strategies to reduce the computational load become especially relevant for higher dimensional settings, such as image classifications, which will be explored in a future work.

\section{Conclusion}
We presented a thorough qualitative analysis of a neural ODE classification task from the perspective of finite-time Lyapunov exponents. As shown for two prototypical vector fields of different structure, these dynamical quantities prove to be convenient and versatile organizers of the input-output dynamics. FTLEs allow a deeper insight into the model's realization of the classification process, locating invariant sets and distinguishing the two phases of early and late dynamics. This knowledge was applied to robustify neural ODEs, by introducing a novel regularization term that manipulates the dynamics by suppressing large FTLEs in the early phase.

\subsection*{Acknowledgments}
CK is supported by a Lichtenberg Professorship. TW is supported by the Austrian Science Fund (FWF) 10.55776/J4681. We thank Elena Queirolo for the insightful discussions and helpful suggestions regarding this work.

\bibliography{25KQW} 
\bibliographystyle{abbrv} 

\end{document}